\pdfoutput=1
\RequirePackage{ifpdf}
\ifpdf 
\documentclass[pdftex]{sigma}
\else
\documentclass{sigma}
\fi

\usepackage{bm}
\usepackage{tikz}
\usetikzlibrary{calc}
\usetikzlibrary{decorations.pathreplacing,calligraphy}

\usepackage{listings}
\lstdefinelanguage{Sage}[]{Python}
{morekeywords={False,sage,True},sensitive=true}
\definecolor{dblackcolor}{rgb}{0.0,0.0,0.0}
\definecolor{dbluecolor}{rgb}{0.01,0.02,0.7}
\definecolor{dgreencolor}{rgb}{0.2,0.4,0.0}
\definecolor{dgraycolor}{rgb}{0.30,0.3,0.30}

\def\zz#1{%
 \ifx\zz#1\else
 #1\linebreak[1]\expandafter\zz
 \fi}

\numberwithin{equation}{section}

\newtheorem{thm}{Theorem}[section]
\newtheorem{prop}[thm]{Proposition}

\newtheorem{lem}[thm]{Lemma}
\newtheorem{cor}[thm]{Corollary}
\newtheorem{conj}[thm]{Conjecture}

\theoremstyle{definition}
\newtheorem{defn}[thm]{Definition}
\newtheorem{eg}[thm]{Example}
\newtheorem{rmk}[thm]{Remark}

\newcommand{\Aa}{\mathbb{A}}
\newcommand{\CC}{\mathbb{C}}

\newcommand{\PP}{\mathbb{P}}

\newcommand{\RR}{\mathbb{R}}
\newcommand{\TT}{\mathbb{T}}
\newcommand{\VV}{\mathbb{V}}
\newcommand{\ZZ}{\mathbb{Z}}

\newcommand{\setupCayleynodes}{
 \draw[gray!30] (0,0) circle (1cm);
 \draw[gray!30] (0,0) circle (0.5176cm);
 \foreach \i in {0,...,11} {
 \node[gray!30] at ({cos(30*\i)},{sin(30*\i)}) {\tiny $\bullet$};
 \node[gray!30] at ({sqrt(2-sqrt(3))*cos(30*\i+15)},{sqrt(2-sqrt(3))*sin(30*\i+15)}) {\tiny $\bullet$};
 \node[gray!30] at (0,0) {\tiny $\bullet$};
 }

 \node[inner sep=-3pt] (a0) at ({cos( 0)},{sin( 0)}) {};
 \node[inner sep=-3pt] (a1) at ({cos( 30)},{sin( 30)}) {};
 \node[inner sep=-3pt] (a2) at ({cos( 60)},{sin( 60)}) {};
 \node[inner sep=-3pt] (a3) at ({cos( 90)},{sin( 90)}) {};
 \node[inner sep=-3pt] (a4) at ({cos(120)},{sin(120)}) {};
 \node[inner sep=-3pt] (a5) at ({cos(150)},{sin(150)}) {};
 \node[inner sep=-3pt] (a6) at ({cos(180)},{sin(180)}) {};
 \node[inner sep=-3pt] (a7) at ({cos(210)},{sin(210)}) {};
 \node[inner sep=-3pt] (a8) at ({cos(240)},{sin(240)}) {};
 \node[inner sep=-3pt] (a9) at ({cos(270)},{sin(270)}) {};
 \node[inner sep=-3pt] (a10) at ({cos(300)},{sin(300)}) {};
 \node[inner sep=-3pt] (a11) at ({cos(330)},{sin(330)}) {};

 \node[inner sep=-3pt] (b0) at ({cos( 0+15)*sqrt(2-sqrt(3))},{sin( 0+15)*sqrt(2-sqrt(3))}) {};
 \node[inner sep=-3pt] (b1) at ({cos( 30+15)*sqrt(2-sqrt(3))},{sin( 30+15)*sqrt(2-sqrt(3))}) {};
 \node[inner sep=-3pt] (b2) at ({cos( 60+15)*sqrt(2-sqrt(3))},{sin( 60+15)*sqrt(2-sqrt(3))}) {};
 \node[inner sep=-3pt] (b3) at ({cos( 90+15)*sqrt(2-sqrt(3))},{sin( 90+15)*sqrt(2-sqrt(3))}) {};
 \node[inner sep=-3pt] (b4) at ({cos(120+15)*sqrt(2-sqrt(3))},{sin(120+15)*sqrt(2-sqrt(3))}) {};
 \node[inner sep=-3pt] (b5) at ({cos(150+15)*sqrt(2-sqrt(3))},{sin(150+15)*sqrt(2-sqrt(3))}) {};
 \node[inner sep=-3pt] (b6) at ({cos(180+15)*sqrt(2-sqrt(3))},{sin(180+15)*sqrt(2-sqrt(3))}) {};
 \node[inner sep=-3pt] (b7) at ({cos(210+15)*sqrt(2-sqrt(3))},{sin(210+15)*sqrt(2-sqrt(3))}) {};
 \node[inner sep=-3pt] (b8) at ({cos(240+15)*sqrt(2-sqrt(3))},{sin(240+15)*sqrt(2-sqrt(3))}) {};
 \node[inner sep=-3pt] (b9) at ({cos(270+15)*sqrt(2-sqrt(3))},{sin(270+15)*sqrt(2-sqrt(3))}) {};
 \node[inner sep=-3pt] (b10) at ({cos(300+15)*sqrt(2-sqrt(3))},{sin(300+15)*sqrt(2-sqrt(3))}) {};
 \node[inner sep=-3pt] (b11) at ({cos(330+15)*sqrt(2-sqrt(3))},{sin(330+15)*sqrt(2-sqrt(3))}) {};}

\newcommand{\setupnodes}{
 \draw[gray!30] (0,0) circle (1cm);
 \draw[gray!30] (0,0) circle (0.653cm);
 \draw[gray!30] (0,0) circle (0.347cm);
 \foreach \i in {0,...,17} {
 \node[gray!30] at ({cos(20*\i)},{sin(20*\i)}) {\tiny $\bullet$};
 \node[gray!30] at ({(1-2*sin(10))*cos(20*\i)},{(1-2*sin(10))*sin(20*\i)}) {\tiny $\bullet$};
 \node[gray!30] at ({0.347*cos(20*\i)},{0.347*sin(20*\i)}) {\tiny $\bullet$};
 \node[gray!30] at (0,0) {$\bullet$};
 }

 \node (a0) at ({cos( 0)},{sin( 0)}) {};
 \node (a1) at ({cos( 20)},{sin( 20)}) {};
 \node (a2) at ({cos( 40)},{sin( 40)}) {};
 \node (a3) at ({cos( 60)},{sin( 60)}) {};
 \node (a4) at ({cos( 80)},{sin( 80)}) {};
 \node (a5) at ({cos(100)},{sin(100)}) {};
 \node (a6) at ({cos(120)},{sin(120)}) {};
 \node (a7) at ({cos(140)},{sin(140)}) {};
 \node (a8) at ({cos(160)},{sin(160)}) {};
 \node (a9) at ({cos(180)},{sin(180)}) {};
 \node (a10) at ({cos(200)},{sin(200)}) {};
 \node (a11) at ({cos(220)},{sin(220)}) {};
 \node (a12) at ({cos(240)},{sin(240)}) {};
 \node (a13) at ({cos(260)},{sin(260)}) {};
 \node (a14) at ({cos(280)},{sin(280)}) {};
 \node (a15) at ({cos(300)},{sin(300)}) {};
 \node (a16) at ({cos(320)},{sin(320)}) {};
 \node (a17) at ({cos(340)},{sin(340)}) {};

 \node (b0) at ({0.653*cos( 0)},{0.653*sin( 0)}) {};
 \node (b1) at ({0.653*cos( 20)},{0.653*sin( 20)}) {};
 \node (b2) at ({0.653*cos( 40)},{0.653*sin( 40)}) {};
 \node (b3) at ({0.653*cos( 60)},{0.653*sin( 60)}) {};
 \node (b4) at ({0.653*cos( 80)},{0.653*sin( 80)}) {};
 \node (b5) at ({0.653*cos(100)},{0.653*sin(100)}) {};
 \node (b6) at ({0.653*cos(120)},{0.653*sin(120)}) {};
 \node (b7) at ({0.653*cos(140)},{0.653*sin(140)}) {};
 \node (b8) at ({0.653*cos(160)},{0.653*sin(160)}) {};
 \node (b9) at ({0.653*cos(180)},{0.653*sin(180)}) {};
 \node (b10) at ({0.653*cos(200)},{0.653*sin(200)}) {};
 \node (b11) at ({0.653*cos(220)},{0.653*sin(220)}) {};
 \node (b12) at ({0.653*cos(240)},{0.653*sin(240)}) {};
 \node (b13) at ({0.653*cos(260)},{0.653*sin(260)}) {};
 \node (b14) at ({0.653*cos(280)},{0.653*sin(280)}) {};
 \node (b15) at ({0.653*cos(300)},{0.653*sin(300)}) {};
 \node (b16) at ({0.653*cos(320)},{0.653*sin(320)}) {};
 \node (b17) at ({0.653*cos(340)},{0.653*sin(340)}) {};

 \node (c0) at ({0.347*cos( 0)},{0.347*sin( 0)}) {};
 \node (c1) at ({0.347*cos( 20)},{0.347*sin( 20)}) {};
 \node (c2) at ({0.347*cos( 40)},{0.347*sin( 40)}) {};
 \node (c3) at ({0.347*cos( 60)},{0.347*sin( 60)}) {};
 \node (c4) at ({0.347*cos( 80)},{0.347*sin( 80)}) {};
 \node (c5) at ({0.347*cos(100)},{0.347*sin(100)}) {};
 \node (c6) at ({0.347*cos(120)},{0.347*sin(120)}) {};
 \node (c7) at ({0.347*cos(140)},{0.347*sin(140)}) {};
 \node (c8) at ({0.347*cos(160)},{0.347*sin(160)}) {};
 \node (c9) at ({0.347*cos(180)},{0.347*sin(180)}) {};
 \node (c10) at ({0.347*cos(200)},{0.347*sin(200)}) {};
 \node (c11) at ({0.347*cos(220)},{0.347*sin(220)}) {};
 \node (c12) at ({0.347*cos(240)},{0.347*sin(240)}) {};
 \node (c13) at ({0.347*cos(260)},{0.347*sin(260)}) {};
 \node (c14) at ({0.347*cos(280)},{0.347*sin(280)}) {};
 \node (c15) at ({0.347*cos(300)},{0.347*sin(300)}) {};
 \node (c16) at ({0.347*cos(320)},{0.347*sin(320)}) {};
 \node (c17) at ({0.347*cos(340)},{0.347*sin(340)}) {};}

\DeclareMathOperator{\Spec}{Spec}

\DeclareMathOperator{\Gr}{Gr}

\DeclareMathOperator{\OGr}{OGr}

\DeclareMathOperator{\Dih}{Dih}

\DeclareMathOperator{\NS}{NS}

\DeclareMathOperator{\dP}{dP}

\DeclareMathOperator{\Nef}{Nef}

\begin{document}

\allowdisplaybreaks

\newcommand{\arXivNumber}{2310.10223}

\renewcommand{\PaperNumber}{033}

\FirstPageHeading

\ShortArticleName{A Laurent Phenomenon for the Cayley Plane}

\ArticleName{A Laurent Phenomenon for the Cayley Plane}

\Author{Oliver DAISEY and Tom DUCAT}

\AuthorNameForHeading{O.~Daisey and T.~Ducat}

\Address{Department of Mathematical Sciences, Durham University, Upper Mountjoy Campus,\\ Stockton Road, Durham DH1 3LE, UK}
\Email{\href{mailto:oliver.j.daisey@durham.ac.uk}{oliver.j.daisey@durham.ac.uk}, \href{mailto:tomaducat@gmail.com}{tomaducat@gmail.com}}
\URLaddress{\url{https://oliverdaisey.github.io/},\newline
\hspace*{10.5mm}\url{https://sites.google.com/site/tomducatmaths}}

\ArticleDates{Received October 22, 2023, in final form April 11, 2024; Published online April 15, 2024}

\Abstract{We describe a Laurent phenomenon for the Cayley plane, which is the homo\-geneous variety associated to the cominuscule representation of $E_6$. The corresponding Laurent phenomenon algebra has finite type and appears in a natural sequence of LPAs indexed by the $E_n$ Dynkin diagrams for $n\leq6$. We conjecture the existence of a further finite type LPA, associated to the Freudenthal variety of type $E_7$.}


\Classification{13F60; 14M17}

\section{Introduction}

\subsection{Background}

In this paper we consider the (complex) Cayley plane, which is the cominuscule homogeneous space \[\mathbb{OP}^2 = E_6/P_6.\] The Cayley plane is a 16-dimensional algebraic variety, with a projective `octonic spinor embedding' $\mathbb{OP}^2\subset \mathbb{P}^{26}$ of codimension 10. We let $\mathcal A_6 := \CC\big[\mathbb{OP}^2\big]$ denote the homogeneous coordinate ring of the Cayley plane with respect to this embedding.

Despite the name, the Cayley plane was first discovered by Ruth Moufang in 1933 \cite{m}. It was named after Cayley since it can be realised as the projective plane over the octonions.

{\bf Laurent phenomenon algebras.}
\emph{Laurent phenomenon algebras} (LPAs) were introduced by Lam and Pylyavskyy \cite{lp} as a generalisation of cluster algebras. As the name suggests, they possess the same Laurent phenomenon property of a cluster algebra but with the flexibility of having much more general exchange polynomials (rather than the more restrictive binomial exchange relations of a cluster algebra).

{\bf Cluster algebra structures on homogeneous spaces.}
A cluster structure (or more generally, a LPA structure) on the coordinate ring of an algebraic variety $\mathcal V$ need not be uniquely determined, since it depends on a choice of anticanonical divisor $\mathcal D\subset \mathcal V$, on which a set of \emph{frozen coefficients} vanish. The simplest examples of cluster structures are \emph{finite type cluster structures}, in which the cluster algebra has only finitely many seeds, although these are rather uncommon since cluster algebras are typically not of finite type.

The existence of cluster algebra structures on the coordinate rings of homogeneous spaces have been an active area of study for some time, beginning with the case Grassmannians~\cite{scott} and followed by the case of partial flag varieties by Gei{\ss}, Leclerc and Schr\"oer~\cite{GLS}. Their construction provides cluster algebra structures for both the Cayley plane and the Freudenthal variety, albeit not ones of finite type.

Given the non-uniqueness of LPA structures in general, it is quite possible (as we show for the Cayley plane in this paper) that a homogeneous variety may admit a finite type LPA structure, even when it only has cluster algebra structures that are not of finite type.

{\bf Mirror symmetry for cluster varieties.}
Cluster varieties have played an important role in the study of mirror symmetry, since they have a well-known mirror construction (the duality between so-called $\mathcal A$- and $\mathcal X$-cluster varieties). In particular, homogeneous spaces have become a~fruitful set of examples for studying various aspects of mirror symmetry, due to the existence of the aforementioned cluster structures. In particular, Spacek and Wang~\cite{Spacek-Wang} recently studied mirrors for both the Cayley plane and the Freudenthal variety using the cluster algebra structures of \cite{GLS}.

{\bf Mirror symmetry for LPAs.}
An LPA is the coordinate ring of a maximal log Calabi--Yau variety and, as such, it is still expected to have a mirror (cf.\ \cite[Conjecture 2.5]{duc}), although an explicit construction for the mirror of an LPA is not yet known. A general approach towards constructing a mirror follows the Gross--Siebert program, which involves building a scattering diagram and then counting broken lines in this scattering diagram to compute coefficients in the equations that define the mirror algebra. In general there may be infinitely many broken lines contributing to these counts, and the mirror can only be defined formally due to convergence issues. However, LPAs of finite type are of particular interest because the corresponding scattering diagram will have only finitely many walls and chambers, and hence only finitely many broken lines contributing to the coefficients of a given equation.

\subsection{Summary of the paper}

{\bf Main result.}
Our main result, Theorem~\ref{thm!main-result}, is the description of a LPA structure of finite type on the ring $\mathcal A_6$. This finite type LPA is of rank 5, and it has 264 seeds and 32 cluster variables.

The key to constructing this LPA structure lies in deriving an initial seed (see Proposition~\ref{prop!initial-seed}) which is compatible with the symmetry of the action of a Coxeter rotation on $\mathcal A_6$. Once we have discovered this seed, the proof of Theorem~\ref{thm!main-result} follows by plugging our seed into the Sage package \textsc{LPASeed} (written by the first author \cite{sagecode}) and verifying the result.

From the output of our code we also verify that our finite type LPA has a positivity phenomenon in Corollary~\ref{cor!positivity}.

{\bf Contents.}
The material in this paper is divided into the following sections.
\begin{itemize}\itemsep=0pt
\item[\S2] A recap of Lam and Pylyavskyy's Laurent phenomenon algebras.
\item[\S3] The description of a sequence of homogeneous varieties $\mathcal V_n$ that we consider in this paper, corresponding to a sequence of type $E_n$ Dynkin diagrams.
\item[\S4] A recap of the finite type LPA structure on the homogeneous coordinate ring of $\mathcal V_n$ for two simpler cases in our sequence. Namely,
\begin{itemize}\itemsep=0pt
\item[{}]{$E_4$ case:} $\mathcal V_4=\Gr(2,5)$ has a finite type cluster algebra structure of type $A_2$, and
\item[{}]{$E_5$ case:} $\mathcal V_5=\OGr(5,10)$ has a finite type LPA structure (considered in \cite{duc}).
\end{itemize}
\item[\S5] Our main $E_6$ case. By generalising the examples of Section~5, we construct a finite type LPA structure on the homogeneous coordinate ring of the Cayley plane $\mathcal V_6=\mathbb{OP}^2$.
\item[\S6] A conjecture and some limited progress on the $E_7$ case.
\item[\S{A}] Examples of the Sage code.
\end{itemize}

\section{Laurent phenomenon algebras}

\subsection{The initial seed}

In the most general setting we consider a coefficient ring $A$, which we take to be an integral domain, and $\mathcal{F}$ a degree $n$ transcendental field extension of $\text{Frac}(A)$.

\begin{defn}
	A \textit{seed} (of rank $n$) consists of a pair $(\bm{x}, \bm{F})$, where
	\begin{enumerate}\itemsep=0pt
		\item[(1)] $\bm{x} = \{x_1, \dots, x_n\}$, the \textit{cluster}, is a transcendence basis for $\mathcal{F}$. The variables $x_1, \dots, x_n$ are called \textit{cluster variables}, and
		\item[(2)] $\bm{F} = \{F_1, \dots, F_n\} \subset A[x_1, \dots, x_n]$ is a collection of $n$ polynomials, called \textit{exchange polynomials}, with the following properties:
		\begin{itemize}\itemsep=0pt
			\item[(LP1)] The exchange polynomial $F_i$ is considered to correspond to the cluster variable $x_i$, and $F_i\in A[x_1, \dots, \widehat{x_i}, \dots, x_n]$ does not depend on $x_i$, and
			\item[(LP2)] each $F_i$ is irreducible and not divisible by any of the cluster variables $x_j$.
		\end{itemize}
	\end{enumerate}
\end{defn}

\begin{defn}\label{def!exchange-laurent-poly}
Given a seed $S = (\bm{x}, \bm{F})$, the \textit{exchange Laurent polynomial} associated to $F_i$ is the Laurent polynomial $\hat{F}_i$ defined according to the following rule: we set $\hat{F}_i = M F_i$ where ${M = x_1^{a_1} \cdots \widehat{x_i} \cdots x_n^{a_n}}$ for some $a_1, \dots , a_{i-1}, a_{i+1}, \dots , a_n \in \mathbb{Z}_{\leq 0}$, and the power $a_j$ of $x_j$ in the denominator of $\hat{F}_i$ is $-m$, where $m$ is the maximal power such that $F_{i} |_{x_j \leftarrow F_j/x_j}$ is divisible by~$F_j^m$.
\end{defn}

The exchange Laurent polynomials are uniquely determined from the exchange polynomials, and vice versa. Indeed, to obtain the exchange Laurent monomials, one sets the power of $x_j$ in the Laurent monomial denominator $M$ equal to the largest power of $F_j$ that divides $F_i$ upon the substitution $x_j \leftarrow F_j / x_j$. To get the exchange polynomials back, one multiplies $\hat{F}_i$ by the unique monomial $M$ (up to a unit) such that $M\hat{F}_i$ is an irreducible polynomial. One important difference to the classical cluster algebra setup is that we use the exchange Laurent polynomials to determine new cluster variables, as opposed to just the exchange polynomials.

\subsection{Mutation procedure}
To mutate a seed $S$ at an index $i$, one proceeds with the following (non-deterministic, cf.\ Remark~\ref{rmk!non-det}) process:
\begin{enumerate}\itemsep=0pt
 \item[(1)] The cluster variables of the mutated seed $\mu_i(S) = ( \mu_i(\bm{x}), \mu_i(\bm{F}) )$ are the same, except we replace $x_i$ according to the exchange relation prescribed by its exchange Laurent polynomial. That is, $\mu_i(\bm{x}) = \{x_1, \dots, x_i', \dots, x_n\}$, where $x_i x_i' = \hat{F}_i$.
 Note that we use the exchange \textit{Laurent} polynomial in this relation.
 \item[(2)] The exchange polynomials $\mu_i(\bm{F}) = \{F_1',\ldots,F_n'\}$ are updated as follows.

 The exchange polynomial $F_{i}'$ corresponding to $x_i'$ remains the same, that is $F_i' := F_i$.

 The exchange polynomials $F_j '$ for $j \neq i$ are defined through the following procedure. If $F_j$ does not depend on $x_i$, then we set $F_j '$ to be any polynomial satisfying $F_j ' = u_j F_j$, where~$u_j$ is a unit in $A$. If $F_j$ does depend on $x_i$, then we perform the following:
 \begin{itemize}\itemsep=0pt
 \item \textit{Substitution step.} We set
 \begin{equation*}
 ( F_j ' )^* = F_j |_{x_i \leftarrow (\hat{F_i} |_{x_{j} \leftarrow 0}) / x_i'}.
 			\end{equation*}
 The substitution is well defined since, by \cite[Lemma~2.7]{lp}, 
 we see that $x_i$ cannot appear in the denominator of $\hat{F}_{j}$.
 \item \textit{Cancellation step.} We divide out by any common factors that $(F_j ')^*$ shares with $\hat{F_i}|_{x_j \leftarrow 0}$. This then defines $F_j '$ up to a monomial multiplier.
 \item \textit{Normalisation step.} We multiply through by a monomial in $x_1, \dots, x_i', \dots, x_n$ to make~$F_j '$ satisfy~(LP1) and~(LP2) as an exchange polynomial in~$S_i$. Such a monomial will be uniquely defined only up to a unit. Thus $F_j '$ is only defined up to a unit multiplier~$u_j$.
 \end{itemize}
 \end{enumerate}

 \begin{eg} \label{eg:mutation}
	Consider the initial seed with exchange variables $\{x_1, x_2\}$ and corresponding exchange polynomials $\{F_1, F_2\} = \{1 + x_2, 1 + x_1\}$. One checks that the exchange Laurent polynomials have trivial denominators, so that $\{F_1, F_2\} = \big\{\hat{F_1}, \hat{F_2}\big\}$. Mutating at $x_1$ replaces $x_1$ with
	\[
	x_3 := \frac{1+x_2}{x_1}
	\]
	and leaves $x_2$ invariant. The exchange polynomial at $x_2$ is changed by this mutation. The substitution step yields
	\[
	(F_2')^* = 1 + \frac{1}{x_3}.
	\]
	There are no common factors to cancel out. Finally, multiplying by $x_3$ gives an exchange polynomial that satisfies (LP1) and (LP2), so we have $F_2' = 1 + x_3$. The mutation is therefore
	\[
	( \{x_1, x_2\}, \{1 + x_2, 1 + x_1\} ) \mapsto \left( \left\{x_3 = \frac{1 + x_2}{x_1}, x_2\right\}, \{1 + x_2, 1 + x_3\} \right).
	\]
\end{eg}

\begin{rmk}\label{rmk!non-det}
The reason that the process is non-deterministic, is the choice of unit multipliers~$u_i$ for the mutated exchange polynomials $F_i'$. If $u_i$ is the unit multiplier for $F_i$ chosen in the mutation~${\mu_i\colon S\to \mu_i(S)}$, we assume that $u_i^{-1}$ is chosen for the mutation $\mu_i\colon \mu_i(S)\to \mu_i(\mu_i(S))$. Thus mutations are involutive, meaning that for any index $i$, the seed $\mu_i(\mu_i(S))$ is equal to $S$.

In practice we assume that all unit multipliers are equal to 1, and identify two seeds whenever they are equivalent, as we now define.
\end{rmk}

\begin{defn}
	Two seeds $S = (\bm{x},\bm{F})$ and $S' = (\bm{x}',\bm{F}')$ of rank $n$ are said to be \textit{equivalent} if for each $1 \leq i \leq n$, there exists units $\mu_i, \tau_i \in \mathcal{F}$ such that $x_i = \mu_i x_i'$ and $F_i = \tau_i F_i'$.
\end{defn}

\subsection{Obtaining the LPA}

In the literature, an LPA is defined as a pair $\left(\mathcal{A}, \left\{S_i\right\}_{i\in I}\right)$, where $\mathcal A$ is a subring of the ambient field $\mathcal F$, and $\left\{ S_i \right\}_{i\in I}$ is a distinguished collection of seeds with cluster variables in $\mathcal F$. The cluster variables generate $\mathcal A$ over $A$, and any two seeds in the collection are mutation-equivalent. For our purposes, it is convenient to present the following definition:
\begin{defn}
	Given a seed $S = (\bm{x}, \bm{F})$, we define the LPA $\mathcal{A}(S)$ generated by $S$ to be the $A$-algebra given by
	\begin{equation*}
		\mathcal{A} = \bigcap_{S_i} A\big[x_{i1}^{\pm1}, \ldots, x_{in}^{\pm1}\big],
	\end{equation*}
	where the index $S_i$ runs over all seeds obtainable through mutation of $S$, and $\{x_{i1}, \dots , x_{in}\}$ is the cluster for $S_i$.
\end{defn}

One can view this definition as the analogue of the \emph{upper cluster algebra}. However, all of the examples we will consider have finitely many cluster variables, and in this case there is no difference between these two definitions for the LPA: both are equal to the ring generated by all of the cluster variables.

\begin{eg}
Continuing with the mutated seed obtained in Example \ref{eg:mutation}, one can compute the mutation at $x_2$ as
\[
 ( \{x_3, x_2\}, \{1 + x_2, 1 + x_3\} ) \mapsto ( \{x_3, x_4\}, \{1 + x_4, 1 + x_3\} ),
\]
where
\[
 x_4 = \frac{1 + x_3}{x_2} = \frac{1 + \frac{1 + x_2}{x_1}}{x_2} = \frac{1 + x_1 + x_2}{x_1x_2}.
\]
Similarly, one obtains by mutation at $x_3$ the seed with cluster variables $\{x_5, x_4\}$ with
\[
 x_5 = \frac{1 + x_4}{x_3} = \frac{x_1(1 + x_2)(1+x_1)}{(1+x_2)x_1x_2} = \frac{1 + x_1}{x_2}
\]
and one can compute that, continuing in this sequence, $x_6 = x_1$, and mutating at $x_5$ returns the initial seed. Thus the corresponding LPA $\mathcal A$ is generated by five cluster variables
\[
 x_1, \qquad x_2, \qquad x_3=\frac{1 + x_2}{x_1}, \qquad x_4=\frac{1 + x_1 + x_2}{x_1x_2}, \qquad x_5=\frac{1 + x_1}{x_2},
\]
and is isomorphic to the ring
\[ \mathcal A \cong \mathbb{C}[x_1,x_2,x_3,x_4,x_5]/I, \]
where $I=(x_5x_2=1+x_1, x_1x_3=1+x_2, \ldots, x_4x_1=1+x_5)$ is the ideal of relations holding between $x_1,\ldots,x_5$.
\end{eg}

{\bf Geometrical interpretation.}
From a more geometrical point of view we consider a LPA~$\mathcal A$ as the ring of regular functions on an affine algebraic variety $U = \Spec \mathcal A$. Each seed $S$ corresponds to the inclusion of a torus chart
\[ S= (\bm x, \bm F), \qquad \TT_S:= \Spec A\big[x_1^{\pm1},\ldots,x_n^{\pm1}\big] \hookrightarrow U, \]
by the Laurent phenomenon $\mathcal A\subseteq A\big[x_1^{\pm1},\ldots,x_n^{\pm1}\big]$. Each mutation then corresponds to a birational map, e.g.,
\[ \mu_1\colon\ \TT_S\dashrightarrow \TT_{\mu_1S}, \qquad \mu_1^*(x_1',x_2,\ldots,x_n) = \big(x_1^{-1}\hat{F}_1,x_2,\ldots,x_n\big), \]
and these torus charts are glued together by identifying points according to these mutations. The LPA (as we have defined it) is then the ring of regular functions on the union of all these seed tori, and $U=\Spec \mathcal A$ is the `affinisation' of $\bigcup_S\TT_S$. This geometrical point of view is described in the context of cluster algebras in \cite[Section~3]{ghk}, and more generally in \cite[Section~2.2]{duc} (see \cite[Remark~2.7]{duc} in particular).

\subsection{Finite type LPAs}
As with the traditional case of cluster algebras, a typical LPA will have infinitely many seeds.

\begin{defn}
If it is only possible to obtain finitely many seeds via mutations of~$S$, then we say that $\mathcal{A}(S)$ is \textit{finite type}, or even that the seed $S$ is finite type.
\end{defn}

One can check that this definition of finite type matches the definition given in the literature (as the finiteness of seeds in the normalisation of $\left(\mathcal{A}, \left\{S_i\right\}_{i\in I} \right)$).

Cluster algebras of finite type have a particularly neat classification, and are in one-to-one correspondence with Dynkin diagrams of finite type. Lam and Pylyavskyy \cite[Theorem~6.6]{lp} classify finite type LPAs of rank 2 and their classification is essentially equivalent to that of cluster algebras of rank 2. However there is no such classification in general, starting with case of LPAs of rank 3. Indeed there are many more finite type LPAs than there are cluster algebras, with Lam and Pylyavskyy showing that the number of finite type LPAs grows exponentially with respect to the rank.

\begin{defn}
	Two seeds $S$, $S'$ are said to be \textit{similar} if there exists a seed $S''$ equivalent to~$S'$, such that $S''$ can be obtained from $S$ by renaming the cluster variables and substituting this renaming into the exchange polynomials. If an LPA has finitely many similarity classes of seeds, we say it is of \textit{finite mutation type}.
\end{defn}

Finite type implies finite mutation type, but not conversely, as shown by Lam and Pylyavskyy in their two-layer brick wall example \cite[Section~7.2]{lp}. One should contrast this situation to the cluster algebra setting in which exceptional quivers arise of infinite type, but finite mutation type.

Borrowing terminology from cluster algebras, we may construct the \textit{exchange graph} of a given seed $S$. The vertices of this graph are given by all seeds obtainable by mutation of $S$, and any two vertices are connected by an edge if one may be obtained from the other (up to similarity) by a single mutation.

\subsection{Description of Sage code}

The first author has implemented the algorithms for computing mutations in Sage. Examples of using the Sage code are available in Appendix~\ref{appendixA}, and a Sage cell for interactive use of the code in a web browser is also available \cite{sagecode}.

{\bf Remarks on the code.}
As in Remark \ref{rmk!non-det}, when we compute a mutation, in the normalisation step, we choose to keep the coefficient of the multiplying monomial equal to 1. Since mutations are defined up to units, this does not change the resulting structure of the exchange graph, neither does it affect the LPA that the cluster variables generate. It does however mean that our statement of the cluster variables is only well-defined up to a unit multiplier.

In our computations and the Sage code, we identify two seeds $S$, $S'$ when they are equivalent.

\section[A type E sequence]{A type $\boldsymbol{E}$ sequence}

The main result of this paper concerns the existence of a finite type LPA, which we view as the~$n=6$ case in a sequence indexed by type $E_n$ Dynkin diagrams.

We consider the Dynkin diagram $E_n$ for $3\leq n\leq 8$ (where $E_3=A_1A_2$, $E_4=A_4$ and $E_5=D_5$) with the nodes labelled according to the convention adopted by Bourbaki:
\begin{center}
\begin{tikzpicture}[scale=1]
 \begin{scope}[xshift=-4cm]
 \node at (0,1) {$E_3$};
 \draw (0,0) -- (1,0);
 \node at (0,0) [label={[label distance=-3pt]below:\scriptsize $1$}] {$\bullet$};
 \node at (1,0) [label={[label distance=-3pt]below:\scriptsize $3$}] {$\bullet$};
 \node at (2,1) [label={[label distance=-3pt]left :\scriptsize $2$}] {$\bullet$};
 \end{scope}
 \begin{scope}
 \node at (0,1) {$E_4$};
 \draw (0,0) -- (2,0) (2,0) -- (2,1);
 \node at (0,0) [label={[label distance=-3pt]below:\scriptsize $1$}] {$\bullet$};
 \node at (1,0) [label={[label distance=-3pt]below:\scriptsize $3$}] {$\bullet$};
 \node at (2,0) [label={[label distance=-3pt]below:\scriptsize $4$}] {$\bullet$};
 \node at (2,1) [label={[label distance=-3pt]left :\scriptsize $2$}] {$\bullet$};
 \end{scope}
 \begin{scope}[xshift=4cm]
 \node at (0,1) {$E_n$};
 \draw (0,0) -- (3.5,0) (4.5,0) -- (5,0) (2,0) -- (2,1);
 \node at (0,0) [label={[label distance=-3pt]below:\scriptsize $1$}] {$\bullet$};
 \node at (1,0) [label={[label distance=-3pt]below:\scriptsize $3$}] {$\bullet$};
 \node at (2,0) [label={[label distance=-3pt]below:\scriptsize $4$}] {$\bullet$};
 \node at (3,0) [label={[label distance=-3pt]below:\scriptsize $5$}] {$\bullet$};
 \node at (4,0) {$\cdots$};
 \node at (5,0) [label={[label distance=-3pt]below:\scriptsize $n$}] {$\bullet$};
 \node at (2,1) [label={[label distance=-3pt]left :\scriptsize $2$}] {$\bullet$};
 \end{scope}
\end{tikzpicture}
\end{center}
We can associate a number of mathematical objects to this sequence:
\begin{enumerate}\itemsep=0pt
\item[(1)] a smooth del Pezzo surface $\dP_n$ of degree $9-n$ \cite[Section~8]{dol}, obtained by blowing up $n$ general points in $\PP^2$,
\item[(2)] the homogeneous space $\mathcal V_n := E_n/P_n$, where $P_n$ is the parabolic subgroup associated to the $n$th node of the Dynkin diagram,
\item[(3)] the $n$-dimensional semiregular Gosset polytope $\Xi_n$ with Coxeter symbol $(n-4)_{21}$ \cite[Section~11.8]{c}.
\end{enumerate}

\subsection{Numerical invariants} \label{sect!numerical}
We collect some numerical invariants associated to this sequence in Table~\ref{table!En}.
\begin{table}[htp]
\centering\renewcommand{\arraystretch}{1.2}
\begin{tabular}{|c|c|c|c|cccc|}\hline
$n$ & $\mathcal V_n$ 			& $\dim \mathcal V_n$ & Coxeter number $h$ 	& $\gamma_{n,1}$& $\gamma_{n,2}$& $\gamma_{n,3}$& $\gamma_{n,2}+\gamma_{n,3}$ \\ \hline
$3$ & $\PP^2\times\PP^1$	& $3$	& $3$		& $6$		& $3$		& $2$ & $5$ \\
$4$ & $\Gr(2,5)$		& $6$	& $5$		& $10$		& $5$		& $5$ & $10$ \\
$5$ & $\OGr(5,10)$		& $10$	& $8$		& $16$		& $10$		& $16$ & $26$ \\
$6$ & $\mathbb{OP}^2$	& $16$	& $12$		& $27$		& $27$		& $72$ & $99$ \\
$7$ & Freudenthal variety	& $27$ & $18$	& $56$		& $126$		& $576$ & $702$ \\ \hline
\end{tabular}
\caption{Numerical invariants associated to the $E_n$ root systems for $3\leq n\leq 7$.}\label{table!En}
\end{table}

According to the three different points of view, for $3\leq n\leq 7$ these numbers can be interpreted in the following ways.
\begin{enumerate}\itemsep=0pt
\item[(1)] The del Pezzo surface $\dP_n$ contains $\gamma_{n,1}$ lines, $\gamma_{n,2}$ conic classes $\xi$ \big(which correspond to conic fibrations $\xi\colon\dP_n\to \PP^1$\big), and $\gamma_{n,3}$ cubic classes $\eta$ \big(which correspond to contractions~${\eta\colon\dP_n\to \PP^2}$\big).
\item[(2)] The homogeneous space $\mathcal V_n\subset\PP^{\gamma_{n,1}-1}$ has an embedding into a projective space of dimension $\gamma_{n,1}-1$ and is cut out by $\gamma_{n,2}$ quadratic equations. It is a Fano variety of the specified dimension and Fano index $h$ (or in other words, $-K_{\mathcal V_n}=\mathcal{O}_{\mathcal V_n}(h)$ for the given embedding).
\item[(3)] The effective cone $\operatorname{Eff}(\dP_n)$ is the cone over $\Xi_n\subset \NS(\dP_n)\cong \RR^{n+1}$, where the polytope~$\Xi_n$ is obtained as the convex hull of the classes of the lines in $\NS(\dP_n)$. It has $\gamma_{n,1}$ vertices and $\gamma_{n,2}+\gamma_{n,3}$ facets, of which $\gamma_{n,2}$ facets are $(n-1)$-dimensional orthoplexes and $\gamma_{n,3}$ facets are $(n-1)$-dimensional simplexes.
\end{enumerate}

\subsection[Coxeter projection of Xi\_n]{Coxeter projection of $\boldsymbol{\Xi_n}$}
\begin{figure}[ht!]
\centering\begin{tikzpicture}[scale=1.5]
 \begin{scope}
 \foreach \i in {1,...,5}
 {
 \draw ({0.382*cos(72*\i+36)},{0.382*sin(72*\i+36)}) -- ({0.382*cos(72*(\i+2)+36)},{0.382*sin(72*(\i+2)+36)});
 \draw ({cos(72*\i)},{sin(72*\i)}) -- ({cos(72*(\i+1))},{sin(72*(\i+1))});
 \draw ({cos(72*\i)},{sin(72*\i)}) -- ({cos(72*(\i+2))},{sin(72*(\i+2))});
 }
 \end{scope}
 \begin{scope}[xshift=2.5cm]
 \foreach \i in {1,...,8}
 {
 \draw ({0.414*cos(45*\i)},{0.414*sin(45*\i)}) -- ({0.414*cos(45*(\i+2))},{0.414*sin(45*(\i+2))});
 \draw ({0.414*cos(45*\i)},{0.414*sin(45*\i)}) -- ({0.414*cos(45*(\i+3))},{0.414*sin(45*(\i+3))});
 \draw ({cos(45*\i)},{sin(45*\i)}) -- ({cos(45*(\i+1))},{sin(45*(\i+1))});
 \draw ({cos(45*\i)},{sin(45*\i)}) -- ({cos(45*(\i+2))},{sin(45*(\i+2))});
 \draw ({cos(45*\i)},{sin(45*\i)}) -- ({cos(45*(\i+3))},{sin(45*(\i+3))});
 \draw ({cos(45*\i)},{sin(45*\i)}) -- ({cos(45*(\i+4))},{sin(45*(\i+4))});
 }
 \end{scope}
 \begin{scope}[xshift=5cm,rotate=15]
\foreach \i in {1,...,12}
{
 \draw ({cos(30*\i+15)},{sin(30*\i+15)}) -- ({cos(30*(\i+1)+15)},{sin(30*(\i+1)+15)});
 \draw ({cos(30*\i+15)},{sin(30*\i+15)}) -- ({cos(30*(\i+2)+15)},{sin(30*(\i+2)+15)});
 \draw ({cos(30*\i+15)},{sin(30*\i+15)}) -- ({cos(30*(\i+3)+15)},{sin(30*(\i+3)+15)});
 \draw ({cos(30*\i+15)},{sin(30*\i+15)}) -- ({cos(30*(\i+4)+15)},{sin(30*(\i+4)+15)});
 \draw ({cos(30*\i+15)},{sin(30*\i+15)}) -- ({cos(30*(\i+5)+15)},{sin(30*(\i+5)+15)});
 \draw ({cos(30*\i+15)},{sin(30*\i+15)}) -- ({cos(30*(\i+6)+15)},{sin(30*(\i+6)+15)});
 \draw ({cos(30*\i)*sqrt(2-sqrt(3))},{sin(30*\i)*sqrt(2-sqrt(3))}) -- ({cos(30*(\i+2))*sqrt(2-sqrt(3))},{sin(30*(\i+2))*sqrt(2-sqrt(3))});
 \draw ({cos(30*\i)*sqrt(2-sqrt(3))},{sin(30*\i)*sqrt(2-sqrt(3))}) -- ({cos(30*(\i+3))*sqrt(2-sqrt(3))},{sin(30*(\i+3))*sqrt(2-sqrt(3))});
 \draw ({cos(30*\i)*sqrt(2-sqrt(3))},{sin(30*\i)*sqrt(2-sqrt(3))}) -- ({cos(30*(\i+5))*sqrt(2-sqrt(3))},{sin(30*(\i+5))*sqrt(2-sqrt(3))});
 \draw ({cos(30*\i)*sqrt(2-sqrt(3))},{sin(30*\i)*sqrt(2-sqrt(3))}) -- ({cos(30*(\i+6))*sqrt(2-sqrt(3))},{sin(30*(\i+6))*sqrt(2-sqrt(3))});
}
 \end{scope}
 \begin{scope}[xshift=7.5cm]
\foreach \i in {1,...,18}
\foreach \j in {1,...,9}
{
 \draw ({cos(20*\i)},{sin(20*\i)}) -- ({cos(20*(\i+\j))},{sin(20*(\i+\j)))});
 \draw ({0.653*cos(20*\i)},{0.653*sin(20*\i)}) -- ({0.653*cos(20*(\i+\j))},{0.653*sin(20*(\i+\j)))});
 \draw ({0.347*cos(20*\i)},{0.347*sin(20*\i)}) -- ({0.347*cos(20*(\i+\j))},{0.347*sin(20*(\i+\j)))});
}
 \end{scope}
\end{tikzpicture}
\caption{Coxeter projections of the polytopes $\Xi_n$ for $n=4,5,6,7$.}\label{fig!polytopes}
\end{figure}
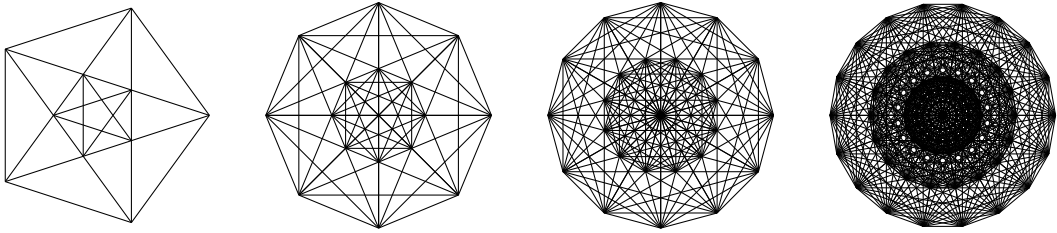

In Figure~\ref{fig!polytopes}, we draw the projection of the polytope $\Xi_n$ onto the Coxeter plane for the cases ${n=4,5,6,7}$. The action of the Coxeter rotation of order $h$ is plainly visible. In these four cases, the vertices of the polytopes are split into orbits of the following sizes:
\[ (10)=2\times(5), \qquad (16)=2\times(8), \qquad (27)=2\times(12)+(3), \qquad (56)=3\times(18)+(2). \]

\subsection{A family of log Calabi--Yau varieties}
Following interpretations (2) and (3) from Section~\ref{sect!numerical}, the $\gamma_{n,1}$ coordinates on the homogeneous space $\mathcal V_n\subset \mathbb{P}^{\gamma_{n,1}-1}$ can be placed in one-to-one correspondence with vertices of $\Xi_n$. Given that the Fano index of $\mathcal V_n$ is equal to the Coxeter number $h$, we can make a natural choice of anticanonical boundary divisor $\mathcal D_n\subset \mathcal V_n$ by taking
$\mathcal D_n = \VV(a_1) + \VV(a_2) + \cdots + \VV(a_h) \in |{-K_{\mathcal V_n}}|$,
where $a_1,\ldots,a_h$ are the coordinates on $\mathcal V_n$ corresponding to the `outside ring' of the Coxeter projection (see Figure~\ref{fig!polytopes}).

\begin{defn}
We let $\mathcal A_n := \CC[\mathcal V_n]$ be the homogeneous coordinate ring of $\mathcal V_n$, with respect to the given embedding, and consider the affine cone $C\mathcal V_n = \Spec \mathcal A_n$. We consider the fibration induced by the projection
\[ \pi\colon\ C\mathcal V_n \to \Aa^{h}_{a_1,\ldots,a_{h}}\]
and let $U_n := \pi^{-1}(\alpha_1, \ldots, \alpha_h)$ denote a general fibre.
\end{defn}

This log Calabi--Yau variety $U_n$ is simply the affine variety obtained by substituting the value~${\alpha_1,\ldots,\alpha_h\in \CC}$ for each coordinate $a_1,\ldots,a_h$ in $\mathcal A_n$ respectively.
This projection $\pi$ spreads out the components of $C\mathcal D_n$ over the coordinate hyperplanes of $\Aa^h$, giving a degenerating family of log Calabi--Yau varieties.

\begin{rmk}\label{rmk!seed-rank}In the language of cluster algebras, the coordinates $a_1,\ldots,a_h$ are \emph{frozen variables} and we will interpret the remaining coordinates as \emph{cluster variables} in a LPA. Moreover, a LPA structure for $\mathcal A_n$ over the base ring $A=\CC[a_1,\ldots,a_h]$ must have rank $r = \dim \mathcal V_n + 1 - h$, since the cluster of each seed will correspond to the inclusion of a torus chart
\[ \Spec A\big[x_1^{\pm1},\ldots,x_r^{\pm1}\big] = \Aa^h_{a_1,\ldots,a_h}\times \big(\CC^\times\big)^r_{x_1,\ldots,x_r} \hookrightarrow C\mathcal V_n \]
that birationally cover $C\mathcal V_n$.
\end{rmk}

\begin{rmk}We do not extend our discussion to include the $E_8$ case, since the numerology of Table~\ref{table!En} breaks down for the homogeneous space $\mathcal V_8=E_8/P_8$. In particular, the Fano index of $\mathcal V_8$ is $29$ (as computed in \cite{snow}), rather than the Coxeter number $h=30$, and thus we do not obtain an anticanonical divisor $\mathcal D_8\subset \mathcal V_8$ in the same way.
\end{rmk}

\section[The Laurent phenomenon for V\_4 and V\_5]{The Laurent phenomenon for $\boldsymbol{\mathcal V_4}$ and $\boldsymbol{\mathcal V_5}$}
\label{section:warm-up-cases}

We briefly summarise the finite type LPA structure on the homogeneous coordinate rings ${\mathcal A_4 = \CC[\mathcal V_4]}$ and $\mathcal A_5 = \CC[\mathcal V_5]$, corresponding to the $E_4$ and $E_5$ cases of our sequence.

\subsection[The E\_4 case]{The $\boldsymbol{E_4}$ case}

The LPA in this case is given by the famous example of the $A_2$ cluster algebra.

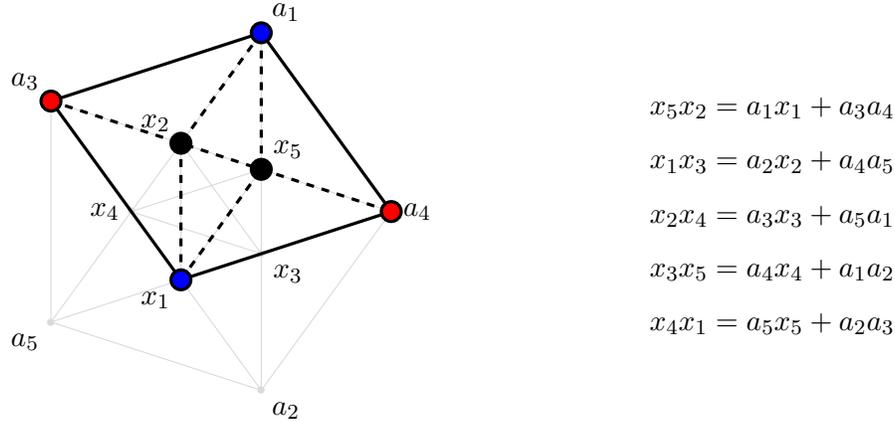
\begin{figure}[ht!]\centering
\begin{tikzpicture}[scale=2.5]
 \node[inner sep=-3pt] (a) at ({cos( 0)},{sin( 0)}) [label={0:$a_4$}] {};
 \node[inner sep=-3pt] (b) at ({cos( 72)},{sin( 72)}) [label={72:$a_1$}] {};
 \node[inner sep=-3pt] (c) at ({cos(144)},{sin(144)}) [label={144:$a_3$}] {};
 \node[inner sep=-3pt] (d) at ({cos(216)},{sin(216)}) [label={216:$a_5$}] {};
 \node[inner sep=-3pt] (e) at ({cos(288)},{sin(288)}) [label={288:$a_2$}] {};

 \node[inner sep=-3pt] (f) at ({0.382*cos( 36)},{0.382*sin( 36)}) [label={36:$x_5$}] {};
 \node[inner sep=-3pt] (g) at ({0.382*cos(108)},{0.382*sin(108)}) [label={108:$x_2$}] {};
 \node[inner sep=-3pt] (h) at ({0.382*cos(180)},{0.382*sin(180)}) [label={180:$x_4$}] {};
 \node[inner sep=-3pt] (i) at ({0.382*cos(252)},{0.382*sin(252)}) [label={252:$x_1$}] {};
 \node[inner sep=-3pt] (j) at ({0.382*cos(324)},{0.382*sin(324)}) [label={324:$x_3$}] {};

 \draw[gray!30] (a) -- (b) -- (c) -- (d) -- (e) -- (a);
 \draw[gray!30] (a) -- (c) -- (e) -- (b) -- (d) -- (a);
 \draw[gray!30] (f) -- (h) -- (j) -- (g) -- (i) -- (f);

 \node[gray!30] (d) at ({cos(216)},{sin(216)}) {\tiny $\bullet$};
 \node[gray!30] (e) at ({cos(288)},{sin(288)}) {\tiny $\bullet$};
 \node[gray!30] (h) at ({0.382*cos(180)},{0.382*sin(180)}) {\tiny $\bullet$};
 \node[gray!30] (j) at ({0.382*cos(324)},{0.382*sin(324)}) {\tiny $\bullet$};

 \draw[very thick] (a) -- (b) -- (c) -- (i) -- (a);
 \draw[very thick,dashed] (b) -- (f) -- (i) -- (g) -- (b);
 \draw[very thick,dashed] (a) -- (c);

 \draw[very thick,fill=red] (a) circle (1.5pt);
 \draw[very thick,fill=red] (c) circle (1.5pt);
 \draw[very thick,fill=blue] (b) circle (1.5pt);

 \draw[very thick,fill=blue] (i) circle (1.5pt);
 \draw[very thick,fill=black] (f) circle (1.5pt);
 \draw[very thick,fill=black] (g) circle (1.5pt);

 \node at (3,0) {\def\arraystretch{1.5} $\begin{array}{l}
 x_5x_2 = a_1x_1 + a_3a_4 \\
 x_1x_3 = a_2x_2 + a_4a_5 \\
 x_2x_4 = a_3x_3 + a_5a_1 \\
 x_3x_5 = a_4x_4 + a_1a_2 \\
 x_4x_1 = a_5x_5 + a_2a_3
 \end{array}$};
\end{tikzpicture}
\caption{The equations of $\mathcal V_4$.}\label{fig!gr25}
\end{figure}

It is convenient to name the frozen variables $a_1,\ldots,a_5$ and the non-frozen variables $x_1,\ldots,x_5$ according to the labelling of $\Xi_4$ shown in Figure~\ref{fig!gr25}. Then the Grassmannian $\mathcal V_4=\Gr(2,5)\subset \PP^9$ is cut out by five quadratic Pl\"ucker equations corresponding to the five octahedral faces of $\Xi_4$.

{\bf Structure of the equations.}
The five quadratic equations have a common structure in that they are comprised of three monomials, each one of which is a product of two opposite vertices in the corresponding octahedron. A coherent choice of signs for these equations is determined by the following \emph{positivity rule}
\begin{gather} x_ix_{i+2}=\text{positive sum of the other monomials}, \label{eq!positivity} \end{gather}
where $x_ix_{i+2}$ is the monomial corresponding to the pair of `internal' vertices of the projected octahedron, and the right hand side comprises of the monomials corresponding to all pairs of `external' vertices.

{\bf Initial seed.}
We let $A= \CC[a_1,\ldots,a_5]$ be the ring generated by the frozen variables and, by Remark~\ref{rmk!seed-rank}, a LPA structure on $\mathcal A_4$ will have rank 2.

As is well-known, each of the following Pl\"ucker coordinates
\[ x_3 = \frac{a_2x_2 + a_4a_5}{x_1}, \qquad x_4 = \frac{a_5a_1x_1 + a_2a_3x_2 + a_3a_4a_5}{x_1x_2}, \qquad x_5 = \frac{a_1x_1 + a_3a_4}{x_2} \]
can be expressed as Laurent polynomials in $A\big[x_1^{\pm1},x_2^{\pm1}\big]$. Moreover, an initial seed for the corresponding LPA structure on $\mathcal A_4$ is given by
\[ S = \begin{cases}
x_1, & a_2x_2 + a_4a_5, \\
x_2, & a_1x_1 + a_3a_4.
\end{cases} \]
Mutating $S$ at $x_1$ gives an almost identical seed (up to reordering) where the only difference is that all the indices of all the variables $x_i$ and $a_i$ have been shifted by $i\mapsto i+1\mod 5$.

{\bf Exchange graph.}
The exchange graph of $\mathcal A_4$ is a pentagon, with vertices labelled by the five possible clusters $\{x_i,x_{i+1}\}$ for all $i\in \ZZ/5\ZZ$ and edges by the five possible mutations~${\{x_{i-1},x_i\}\to\{x_i,x_{i+1}\}}$.

{\bf Positivity.}
The ring $\mathcal A_4$ also has a curious property known as the \emph{positivity phenomenon}: every coefficient in the Laurent expansion of every cluster variable is positive, as well as every coefficient in the exchange polynomials of every seed.

\subsection[The E\_5 case]{The $\boldsymbol{E_5}$ case}
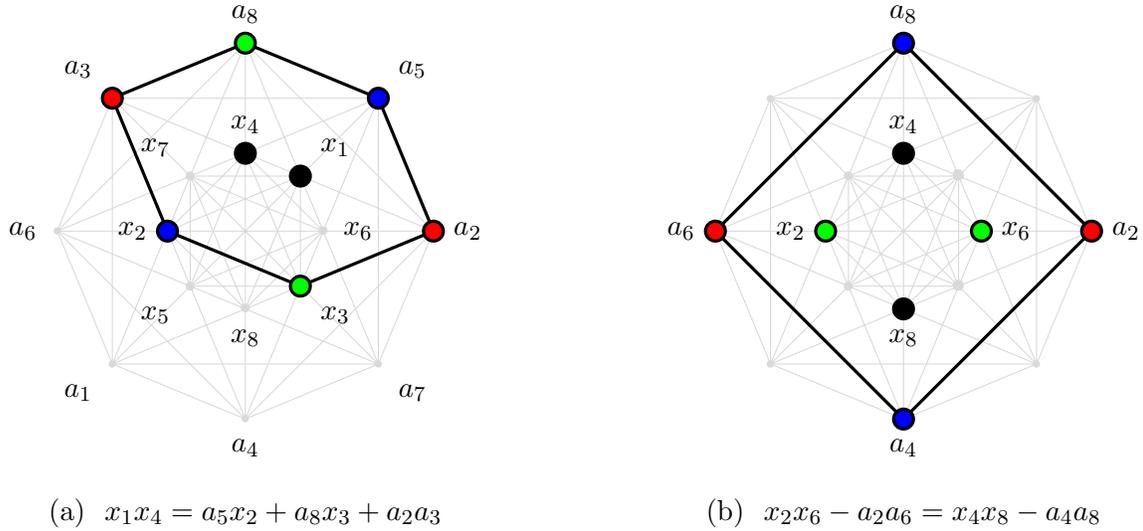
\begin{figure}[ht!]\centering
\begin{tikzpicture}[scale=2.5]
 \node[inner sep=-3pt,gray!30] (a) at ({cos( 0)},{sin( 0)}) {\tiny $\bullet$};
 \node[inner sep=-3pt,gray!30] (b) at ({cos( 45)},{sin( 45)}) {\tiny $\bullet$};
 \node[inner sep=-3pt,gray!30] (c) at ({cos( 90)},{sin( 90)}) {\tiny $\bullet$};
 \node[inner sep=-3pt,gray!30] (d) at ({cos(135)},{sin(135)}) {\tiny $\bullet$};
 \node[inner sep=-3pt,gray!30] (e) at ({cos(180)},{sin(180)}) {\tiny $\bullet$};
 \node[inner sep=-3pt,gray!30] (f) at ({cos(225)},{sin(225)}) {\tiny $\bullet$};
 \node[inner sep=-3pt,gray!30] (g) at ({cos(270)},{sin(270)}) {\tiny $\bullet$};
 \node[inner sep=-3pt,gray!30] (h) at ({cos(315)},{sin(315)}) {\tiny $\bullet$};

 \node[inner sep=-3pt,gray!30] (a2) at ({0.414*cos( 0)},{0.414*sin( 0)}) {\scriptsize $\bullet$};
 \node[inner sep=-3pt,gray!30] (b2) at ({0.414*cos( 45)},{0.414*sin( 45)}) {$\bullet$};
 \node[inner sep=-3pt,gray!30] (c2) at ({0.414*cos( 90)},{0.414*sin( 90)}) {$\bullet$};
 \node[inner sep=-3pt,gray!30] (d2) at ({0.414*cos(135)},{0.414*sin(135)}) {\scriptsize $\bullet$};
 \node[inner sep=-3pt,gray!30] (e2) at ({0.414*cos(180)},{0.414*sin(180)}) {$\bullet$};
 \node[inner sep=-3pt,gray!30] (f2) at ({0.414*cos(225)},{0.414*sin(225)}) {\scriptsize $\bullet$};
 \node[inner sep=-3pt,gray!30] (g2) at ({0.414*cos(270)},{0.414*sin(270)}) {\scriptsize $\bullet$};
 \node[inner sep=-3pt,gray!30] (h2) at ({0.414*cos(315)},{0.414*sin(315)}) {$\bullet$};

 \draw[gray!30] (a) -- (b) -- (c) -- (d) -- (e) -- (f) -- (g) -- (h) -- (a);
 \draw[gray!30] (a) -- (c) -- (e) -- (g) -- (a) (b) -- (d) -- (f) -- (h) -- (b);
 \draw[gray!30] (a2) -- (c2) -- (e2) -- (g2) -- (a2) (b2) -- (d2) -- (f2) -- (h2) -- (b2);
 \draw[gray!30] (a) -- (d) -- (g) -- (b) -- (e) -- (h) -- (c) -- (f) -- (a);
 \draw[gray!30] (a2) -- (d2) -- (g2) -- (b2) -- (e2) -- (h2) -- (c2) -- (f2) -- (a2);
 \draw[gray!30] (a) -- (e) (b) -- (f) (c) -- (g) (d) -- (h);

 \draw[very thick] (a) -- (b) -- (c) -- (d) -- (e2) -- (h2) -- (a);

 \node at (a) [label={ 0:$a_2$}] {};
 \node at (b) [label={ 45:$a_5$}] {};
 \node at (c) [label={ 90:$a_8$}] {};
 \node at (d) [label={135:$a_3$}] {};
 \node at (e) [label={180:$a_6$}] {};
 \node at (f) [label={225:$a_1$}] {};
 \node at (g) [label={270:$a_4$}] {};
 \node at (h) [label={-45:$a_7$}] {};
 \node at (a2) [label={ 0:$x_6$}] {};
 \node at (b2) [label={ 45:$x_1$}] {};
 \node at (c2) [label={ 90:$x_4$}] {};
 \node at (d2) [label={135:$x_7$}] {};
 \node at (e2) [label={180:$x_2$}] {};
 \node at (f2) [label={225:$x_5$}] {};
 \node at (g2) [label={270:$x_8$}] {};
 \node at (h2) [label={-45:$x_3$}] {};

 \draw[very thick,fill=red] (a) circle (1.5pt);
 \draw[very thick,fill=red] (d) circle (1.5pt);
 \draw[very thick,fill=blue] (b) circle (1.5pt);
 \draw[very thick,fill=blue] (e2) circle (1.5pt);
 \draw[very thick,fill=green] (c) circle (1.5pt);
 \draw[very thick,fill=green] (h2) circle (1.5pt);
 \draw[very thick,fill=black] (b2) circle (1.5pt);
 \draw[very thick,fill=black] (c2) circle (1.5pt);

 \begin{scope}[xshift=3.5cm]
 \node[inner sep=-3pt,gray!30] (a) at ({cos( 0)},{sin( 0)}) {\tiny $\bullet$};
 \node[inner sep=-3pt,gray!30] (b) at ({cos( 45)},{sin( 45)}) {\tiny $\bullet$};
 \node[inner sep=-3pt,gray!30] (c) at ({cos( 90)},{sin( 90)}) {\tiny $\bullet$};
 \node[inner sep=-3pt,gray!30] (d) at ({cos(135)},{sin(135)}) {\tiny $\bullet$};
 \node[inner sep=-3pt,gray!30] (e) at ({cos(180)},{sin(180)}) {\tiny $\bullet$};
 \node[inner sep=-3pt,gray!30] (f) at ({cos(225)},{sin(225)}) {\tiny $\bullet$};
 \node[inner sep=-3pt,gray!30] (g) at ({cos(270)},{sin(270)}) {\tiny $\bullet$};
 \node[inner sep=-3pt,gray!30] (h) at ({cos(315)},{sin(315)}) {\tiny $\bullet$};

 \node[inner sep=-3pt,gray!30] (a2) at ({0.414*cos( 0)},{0.414*sin( 0)}) {\scriptsize $\bullet$};
 \node[inner sep=-3pt,gray!30] (b2) at ({0.414*cos( 45)},{0.414*sin( 45)}) {$\bullet$};
 \node[inner sep=-3pt,gray!30] (c2) at ({0.414*cos( 90)},{0.414*sin( 90)}) {$\bullet$};
 \node[inner sep=-3pt,gray!30] (d2) at ({0.414*cos(135)},{0.414*sin(135)}) {\scriptsize $\bullet$};
 \node[inner sep=-3pt,gray!30] (e2) at ({0.414*cos(180)},{0.414*sin(180)}) {$\bullet$};
 \node[inner sep=-3pt,gray!30] (f2) at ({0.414*cos(225)},{0.414*sin(225)}) {\scriptsize $\bullet$};
 \node[inner sep=-3pt,gray!30] (g2) at ({0.414*cos(270)},{0.414*sin(270)}) {\scriptsize $\bullet$};
 \node[inner sep=-3pt,gray!30] (h2) at ({0.414*cos(315)},{0.414*sin(315)}) {$\bullet$};

 \draw[gray!30] (a) -- (b) -- (c) -- (d) -- (e) -- (f) -- (g) -- (h) -- (a);
 \draw[gray!30] (a) -- (c) -- (e) -- (g) -- (a) (b) -- (d) -- (f) -- (h) -- (b);
 \draw[gray!30] (a2) -- (c2) -- (e2) -- (g2) -- (a2) (b2) -- (d2) -- (f2) -- (h2) -- (b2);
 \draw[gray!30] (a) -- (d) -- (g) -- (b) -- (e) -- (h) -- (c) -- (f) -- (a);
 \draw[gray!30] (a2) -- (d2) -- (g2) -- (b2) -- (e2) -- (h2) -- (c2) -- (f2) -- (a2);
 \draw[gray!30] (a) -- (e) (b) -- (f) (c) -- (g) (d) -- (h);

 \draw[very thick] (a) -- (c) -- (e) -- (g) -- (a);

 \node at (a) [label={ 0:$a_2$}] {};
 \node at (e) [label={180:$a_6$}] {};
 \node at (c) [label={ 90:$a_8$}] {};
 \node at (g) [label={270:$a_4$}] {};
 \node at (a2) [label={ 0:$x_6$}] {};
 \node at (e2) [label={180:$x_2$}] {};
 \node at (c2) [label={ 90:$x_4$}] {};
 \node at (g2) [label={270:$x_8$}] {};

 \draw[very thick,fill=red] (a) circle (1.5pt);
 \draw[very thick,fill=red] (e) circle (1.5pt);
 \draw[very thick,fill=blue] (c) circle (1.5pt);
 \draw[very thick,fill=blue] (g) circle (1.5pt);
 \draw[very thick,fill=green] (a2) circle (1.5pt);
 \draw[very thick,fill=green] (e2) circle (1.5pt);
 \draw[very thick,fill=black] (c2) circle (1.5pt);
 \draw[very thick,fill=black] (g2) circle (1.5pt);
 \end{scope}

 \node at (0,-1.5) {(a) \def\arraystretch{1.5} $
 x_1x_4 = a_5x_2 + a_8x_3 + a_2a_3
 $};

 \node at (3.5,-1.5) {(b) \def\arraystretch{1.5} $
 x_2x_6 - a_2a_6 = x_4x_8 - a_4a_8
 $};
\end{tikzpicture}
\caption{The equations of $\mathcal V_5$.}\label{fig!ogr5-10}
\end{figure}

This is the LPA studied in \cite{duc}. In this case we label the variables $a_1,\ldots,a_8$ and $x_1,\ldots,x_8$ with~${i\in \ZZ/8\ZZ}$, as in Figure~\ref{fig!ogr5-10}. The orthogonal Grassmannian $\mathcal V_5=\OGr(5,10)\subset \PP^{15}$ is cut out by ten quadratic equations which correspond to the ten octahedral faces of $\Xi_5$. However this time the equations of $\mathcal V_5$ split into one orbit (a) of size eight and one orbit (b) of size two.

{\bf Structure of the equations.}
We can make a coherent choice of minus signs in the equations by asking that the eight (a) equations obey the analogous positivity rule to equation \eqref{eq!positivity}. Doing that uniquely determines the signs in the remaining two (b) equations.

{\bf Initial seed.}
By Remark~\ref{rmk!seed-rank}, a LPA structure on $\mathcal A_5$ will have rank 3. Beginning with $\{x_1,x_2,x_3\}$ as a candidate for an initial cluster, we can check that each of the other $x_i$ can be written as a Laurent polynomial in $x_1$, $x_2$, $x_3$. Thus we might hope that this initial cluster can be used to get an LPA structure on $\mathcal V_5$ which is analogous to the LPA structure on $\mathcal V_4$.

To promote this cluster into a seed, we have to specify what the exchange polynomial $F_i$ corresponding to $x_i$ should be for $i=1,2,3$. The two equations $x_1x_4=\cdots$ and $x_3x_8=\cdots$ give easy and obvious candidates for the exchange polynomials $F_1$ and $F_3$:
\[ F_1 = a_5x_2 + a_8x_3 + a_2a_3, \qquad F_3 = a_4x_1 + a_7x_2 + a_1a_2. \]
However, it is not immediately clear how to write down the exchange polynomial $F_2$.

Since we would like mutation in the LPA to be compatible with the Coxeter symmetry (as it was in the previous case), we can easily work out what $F_2$ should be by considering the mutation $\mu_1\colon \{x_1,x_2,x_3\}\mapsto \{x_2,x_3,x_4\}$, writing down the exchange polynomial $\mu_1F_2=a_6x_3 + a_1x_4 + a_3a_4$ that we expect to see for $x_2$ with respect to this seed, and then mutating back to get $F_2 = \mu_1^{-1}(\mu_1F_2)$.

As seen in \cite{duc}, this gives an initial seed
\[ S = \begin{cases}
x_1, & a_5x_2 + a_8x_3 + a_2a_3, \\
x_2, & a_6x_1x_3 + a_3a_4x_1 + a_8a_1x_3 + a_1a_2a_3, \\
x_3, & a_4x_1 + a_7x_2 + a_1a_2,
\end{cases} \]
and, incredibly, the mutation of this LPA seed is compatible with the $\Dih_8$-symmetry, in the sense that mutating $S$ at $x_1$ returns an identical seed (up to reordering) with all indices shifted by $i\mapsto i+1\mod 8$.

Moreover, mutating $x_2$ gives a quantity $q_1 = x_2^{-1}F_2$ which can be expressed as a quadratic $q_1 = x_1x_5 - a_1a_5 = x_3x_7 - a_3a_7$ in the other variables. If we also let $q_2 = x_2x_6 - a_2a_6 = x_4x_8 - a_4a_8$, then we get a finite type LPA structure with sixteen clusters:
\begin{gather*} \{x_1,x_2,x_3\}, \{x_2,x_3,x_4\}, \{x_3,x_4,x_5\}, \dots, \{x_7,x_8,x_1\}, \\
 \{x_1,x_3,q_1\}, \{x_2,x_4,q_2\}, \{x_3,x_5,q_1\}, \dots, \{x_8,x_2,q_2\}. \end{gather*}

{\bf Exchange graph.}
The exchange graph is the 1-skeleton of a 3-dimensional polytope with sixteen vertices, eight pentagonal faces (corresponding to $x_1,\ldots,x_8$) and two square faces (corresponding to $q_1$, $q_2$). The exchange graph, shown in Figure~\ref{fig:exchange-graph}, can be produced using the Sage code found in Appendix~\ref{appendixA}.

\begin{figure}[th!]
\centering
\begin{tabular}{ c p{2cm} }
 \begin{minipage}{.45\textwidth}
 \includegraphics[scale=0.5]{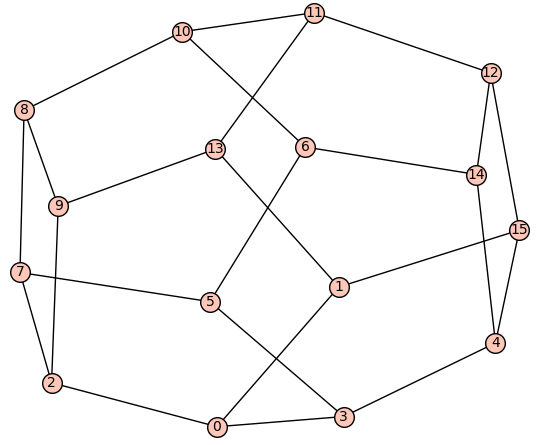}
 \end{minipage}
 &
 $f = \begin{pmatrix}16 \\ 24 \\ 10 \\ 1\end{pmatrix}$
\end{tabular}

\caption{The exchange graph for the LPA in the $\mathcal V_5$ case with f-vector $f$.}\label{fig:exchange-graph}
\end{figure}

{\bf Positivity.}
As with the previous case, $\mathcal A_5$ also has the positivity phenomenon. By an explicit calculation one can check that every Laurent expansion of a cluster variable, as well as every exchange polynomial in every seed, has positive coefficients.\footnote{The reader may be a little bit disturbed by the fact that the two equations of type (b) appear to have negative coefficients. However, after introducing the two new cluster variables $q_1$, $q_2$, they can be rewritten as exchange relations with positive coefficients, e.g., $x_1x_5-a_1a_5 = x_3x_7 - a_3a_7 \implies x_1x_5 = q_1 + a_1a_5$.}

\section[The Laurent phenomenon for V\_6]{The Laurent phenomenon for $\boldsymbol{\mathcal V_6}$}

To describe the projective embedding of the Cayley plane, we must first understand the equations of the projective embedding $\mathcal V_6\subset \PP^{26}$. We do this by thinking of the 27-dimensional representation of $E_6$ in terms of the 27 lines on a cubic surface.

\subsection{The equations of the Cayley plane} \label{sec!cayley-eqns}

We fix a birational model $\pi\colon \dP_6\to \PP^2$ for the smooth cubic surface $\dP_6\subset \PP^3$ obtained as the blowup of six points $p_1,\ldots,p_6\in\PP^2$, and we name the 27 lines in $\dP_6$ according to the following conventions:
\begin{enumerate}\itemsep=0pt
\item[(1)] let $e_i$ be the line corresponding to the exceptional divisor over $p_i$,
\item[(2)] let $\ell_{ij}$ be the line corresponding to the line through $p_i$, $p_j$ and
\item[(3)] let $c_i$ be the line corresponding to the conic through the five points other than $p_i$.
\end{enumerate}

{\bf The polytope $\boldsymbol{\Xi_6}$.}
As mentioned in Section~\ref{sect!numerical}, the polytope $\Xi_6\subset \NS(\dP_6)$ has 27 vertices, corresponding to the 27 lines of $\dP_6$, and $99=27+72$ facets which are one of two types:
\begin{enumerate}\itemsep=0pt
\item[(1)] There are 27 5-dimensional orthoplex facets, corresponding to the 27 extremal rays $.\xi\in \Nef(\dP_6)$ that define conic fibrations $\pi_\xi\colon \dP_6\to \PP^1$. Each face has ten vertices $\ell_1,\ldots,\ell_5$, $\ell_1',\ldots,\ell_5'$ appearing in five `opposite' pairs $\ell_i$, $\ell_i'$ such that $\xi\sim \ell_i + \ell_i'$ for $i=1,\ldots,5$.
\item[(2)] There are 72 5-dimensional simplex facets, corresponding to the 72 extremal rays $\eta\in \Nef(\dP_6)$ that define contractions $\pi_\eta\colon \dP_6\to \PP^2$. Each face has six vertices $\ell_1,\ldots,\ell_6$; the six lines that are contracted by $\pi_\eta$.
\end{enumerate}

\begin{figure}[ht!]\centering
\begin{tikzpicture}[scale=2.8,font=\small]

\foreach \i in {1,...,12}
{
 \draw ({cos(30*\i+15)},{sin(30*\i+15)}) -- ({cos(30*(\i+1)+15)},{sin(30*(\i+1)+15)});
 \draw ({cos(30*\i+15)},{sin(30*\i+15)}) -- ({cos(30*(\i+2)+15)},{sin(30*(\i+2)+15)});
 \draw ({cos(30*\i+15)},{sin(30*\i+15)}) -- ({cos(30*(\i+3)+15)},{sin(30*(\i+3)+15)});
 \draw ({cos(30*\i+15)},{sin(30*\i+15)}) -- ({cos(30*(\i+4)+15)},{sin(30*(\i+4)+15)});
 \draw ({cos(30*\i+15)},{sin(30*\i+15)}) -- ({cos(30*(\i+5)+15)},{sin(30*(\i+5)+15)});
 \draw ({cos(30*\i+15)},{sin(30*\i+15)}) -- ({cos(30*(\i+6)+15)},{sin(30*(\i+6)+15)});
 \draw ({cos(30*\i)*sqrt(2-sqrt(3))},{sin(30*\i)*sqrt(2-sqrt(3))}) -- ({cos(30*(\i+2))*sqrt(2-sqrt(3))},{sin(30*(\i+2))*sqrt(2-sqrt(3))});
 \draw ({cos(30*\i)*sqrt(2-sqrt(3))},{sin(30*\i)*sqrt(2-sqrt(3))}) -- ({cos(30*(\i+3))*sqrt(2-sqrt(3))},{sin(30*(\i+3))*sqrt(2-sqrt(3))});
 \draw ({cos(30*\i)*sqrt(2-sqrt(3))},{sin(30*\i)*sqrt(2-sqrt(3))}) -- ({cos(30*(\i+4))*sqrt(2-sqrt(3))},{sin(30*(\i+4))*sqrt(2-sqrt(3))});
 \draw ({cos(30*\i)*sqrt(2-sqrt(3))},{sin(30*\i)*sqrt(2-sqrt(3))}) -- ({cos(30*(\i+5))*sqrt(2-sqrt(3))},{sin(30*(\i+5))*sqrt(2-sqrt(3))});
 \draw ({cos(30*\i)*sqrt(2-sqrt(3))},{sin(30*\i)*sqrt(2-sqrt(3))}) -- ({cos(30*(\i+6))*sqrt(2-sqrt(3))},{sin(30*(\i+6))*sqrt(2-sqrt(3))});
}
	\draw[fill=white,thick] ({cos(30*0+15)},{sin(30*0+15)}) circle (3pt) node {$a_1$};
	\draw[fill=white,thick] ({cos(30*1+15)},{sin(30*1+15)}) circle (3pt) node {$a_6$};
	\draw[fill=white,thick] ({cos(30*2+15)},{sin(30*2+15)}) circle (3pt) node {$a_{11}$};
	\draw[fill=white,thick] ({cos(30*3+15)},{sin(30*3+15)}) circle (3pt) node {$a_4$};
	\draw[fill=white,thick] ({cos(30*4+15)},{sin(30*4+15)}) circle (3pt) node {$a_9$};
	\draw[fill=white,thick] ({cos(30*5+15)},{sin(30*5+15)}) circle (3pt) node {$a_2$};
	\draw[fill=white,thick] ({cos(30*6+15)},{sin(30*6+15)}) circle (3pt) node {$a_7$};
	\draw[fill=white,thick] ({cos(30*7+15)},{sin(30*7+15)}) circle (3pt) node {$a_{12}$};
	\draw[fill=white,thick] ({cos(30*8+15)},{sin(30*8+15)}) circle (3pt) node {$a_5$};
	\draw[fill=white,thick] ({cos(30*9+15)},{sin(30*9+15)}) circle (3pt) node {$a_{10}$};
	\draw[fill=white,thick] ({cos(30*10+15)},{sin(30*10+15)}) circle (3pt) node {$a_3$};
	\draw[fill=white,thick] ({cos(30*11+15)},{sin(30*11+15)}) circle (3pt) node {$a_8$};
	
	\draw[fill=white,thick] ({cos(30*0)*sqrt(2-sqrt(3))},{sin(30*0)*sqrt(2-sqrt(3))}) circle (3pt) node {$x_1$};
	\draw[fill=white,thick] ({cos(30*1)*sqrt(2-sqrt(3))},{sin(30*1)*sqrt(2-sqrt(3))}) circle (3pt) node {$x_6$};
	\draw[fill=white,thick] ({cos(30*2)*sqrt(2-sqrt(3))},{sin(30*2)*sqrt(2-sqrt(3))}) circle (3pt) node {$x_{11}$};
	\draw[fill=white,thick] ({cos(30*3)*sqrt(2-sqrt(3))},{sin(30*3)*sqrt(2-sqrt(3))}) circle (3pt) node {$x_4$};
	\draw[fill=white,thick] ({cos(30*4)*sqrt(2-sqrt(3))},{sin(30*4)*sqrt(2-sqrt(3))}) circle (3pt) node {$x_9$};
	\draw[fill=white,thick] ({cos(30*5)*sqrt(2-sqrt(3))},{sin(30*5)*sqrt(2-sqrt(3))}) circle (3pt) node {$x_2$};
	\draw[fill=white,thick] ({cos(30*6)*sqrt(2-sqrt(3))},{sin(30*6)*sqrt(2-sqrt(3))}) circle (3pt) node {$x_7$};
	\draw[fill=white,thick] ({cos(30*7)*sqrt(2-sqrt(3))},{sin(30*7)*sqrt(2-sqrt(3))}) circle (3pt) node {$x_{12}$};
	\draw[fill=white,thick] ({cos(30*8)*sqrt(2-sqrt(3))},{sin(30*8)*sqrt(2-sqrt(3))}) circle (3pt) node {$x_5$};
	\draw[fill=white,thick] ({cos(30*9)*sqrt(2-sqrt(3))},{sin(30*9)*sqrt(2-sqrt(3))}) circle (3pt) node {$x_{10}$};
	\draw[fill=white,thick] ({cos(30*10)*sqrt(2-sqrt(3))},{sin(30*10)*sqrt(2-sqrt(3))}) circle (3pt) node {$x_3$};
	\draw[fill=white,thick] ({cos(30*11)*sqrt(2-sqrt(3))},{sin(30*11)*sqrt(2-sqrt(3))}) circle (3pt) node {$x_8$};
	
	\draw[fill=white,thick] (0,0) circle (3pt) node {$z_i$};
	

 \begin{scope}[scale=0.5, xshift = 4.5cm, yshift = 1.25cm]
 \node at (-1,1) {(a)};
 \begin{scope}[rotate=-60]

 \setupCayleynodes

 \draw[thick] (a0) -- (a1) -- (a2) -- (a3) -- (a4) -- (b4) -- (0,0) -- (b11) -- (a0);

 \node[red] at (a0) {$\bullet$};
 \node[red] at (a4) {$\bullet$};
 \node[blue] at (a1) {$\bullet$};
 \node[blue] at (b4) {$\bullet$};
 \node[green] at (a2) {$\bullet$};
 \node[green] at (0,0) {$\bullet$};
 \node[purple] at (a3) {$\bullet$};
 \node[purple] at (b11) {$\bullet$};
 \node at (b1) {$\bullet$};
 \node at (b2) {$\bullet$};
 \end{scope}
 \end{scope}

 \begin{scope}[scale=0.5, xshift = 7cm, yshift = 0cm]
 \node at (-1,1) {(b)};
 \begin{scope}[rotate=45]
 \setupCayleynodes

 \draw[thick] (a8) -- (a10) -- (a11) -- (a1) -- (b3) -- (b5) -- (a8);

 \node at (b0) {$\bullet$};
 \node at (b8) {$\bullet$};
 \node[green] at (b10) {$\bullet$};
 \node[green] at (0,0) {$\bullet$};
 \node[purple] at (b5) {$\bullet$};
 \node[purple] at (a11) {$\bullet$};
 \node[blue] at (b3) {$\bullet$};
 \node[blue] at (a10) {$\bullet$};
 \node[red] at (a1) {$\bullet$};
 \node[red] at (a8) {$\bullet$};
 \end{scope}
 \end{scope}

 \begin{scope}[scale=0.5, xshift = 4.5cm, yshift = -1.25cm]
 \node at (-1,1) {(c)};
 \begin{scope}[rotate=15]
 \setupCayleynodes

 \draw[thick] (a1) -- (a4) -- (a7) -- (a10) -- (a1);

 \node[blue] at (a1) {$\bullet$};
 \node[blue] at (a7) {$\bullet$};
 \node at (0,0) {$\bullet$};
 \node[purple] at (a4) {$\bullet$};
 \node[purple] at (a10) {$\bullet$};
 \node[green] at (b11) {$\bullet$};
 \node[green] at (b5) {$\bullet$};
 \node[red] at (b2) {$\bullet$};
 \node[red] at (b8) {$\bullet$};
 \end{scope}
 \end{scope}
\end{tikzpicture}
\caption{(i) A labelling of the vertices of $\Xi_6$, and (ii) the three types of octahedral face.}\label{fig!cayley-plane}
\end{figure}
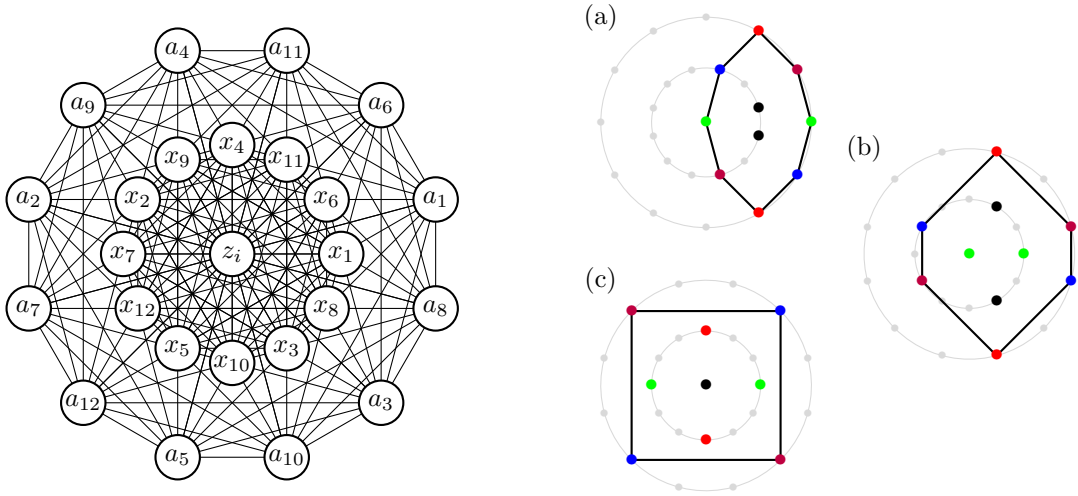

{\bf Action of a Coxeter rotation.}
The 72 roots of the $E_6$ root system in $\NS(\dP_6)$ are given by all possible differences $\ell_i-\ell_j$, where $\ell_i,\ell_j$ are a pair of non-intersecting lines in $\dP_6$. Each root $r_i$ specifies a reflection \[\rho_i(r_j) = r_j + (r_i\cdot r_j)r_i,\] and a Coxeter rotation is an element~${\sigma=\rho_1\rho_2\rho_3\rho_4\rho_5\rho_6}$ of order 12, obtained as a product of the reflections over a set of simple roots. For example, if~$\sigma$ is the Coxeter rotation obtained from the following choice of simple roots
\begin{center}
\begin{tikzpicture}[xscale=2.5]
 \draw (0,0) -- (4,0) (2,0) -- (2,1);
 \node at (0,0) [label={[label distance=-3pt]below:\small $r_1=e_1-e_2$}] {$\bullet$};
 \node at (1,0) [label={[label distance=-3pt]below:\small $r_3=e_2-e_3$}] {$\bullet$};
 \node at (2,0) [label={[label distance=-3pt]below:\small $r_4=e_3-e_4$}] {$\bullet$};
 \node at (3,0) [label={[label distance=-3pt]below:\small $r_5=e_4-e_5$}] {$\bullet$};
 \node at (4,0) [label={[label distance=-3pt]below:\small $r_6=e_5-e_6$}] {$\bullet$};
 \node at (2,1) [label={[label distance=-3pt]left :\small $r_2=\ell_{12}-e_3$}] {$\bullet$};
\end{tikzpicture}
\end{center}
then the action of $\sigma$ on the set of the 27 lines has (ordered) orbits of length 12, 12, 3, as shown in the rows of Table~\ref{table!e6-orbits}. The 27 lines correspond to 27 spinor coordinates of $\mathcal V_6\subset\PP^{26}$, and we name these coordinates $a_i$, $x_j$, $z_k$ for $i,j\in\ZZ/12\ZZ$ and $k\in\ZZ/3\ZZ$ according to which of these three orbits they belong to, as in Table~\ref{table!e6-orbits}.
\begin{table}[htbp]\renewcommand{\arraystretch}{1.2}
$$\begin{array}{|c|cccccc cccccc|} \hline
i & 1 & 2 & 3 & 4 & 5 & 6 & 7 & 8 & 9 & 10 & 11 & 12 \\ \hline
a_i & e_3 & \ell_{34} & \ell_{56} & e_6 & c_2 & e_4 & c_6 & \ell_{16} & \ell_{23} & c_3 & e_5 & c_1 \\
x_i & \ell_{26} & \ell_{24} & \ell_{46} & e_1 & \ell_{45} & e_2 & \ell_{35} & \ell_{15} & \ell_{13} & c_4 & \ell_{12} & c_5 \\
z_i & \ell_{14} & \ell_{36} & \ell_{25} & & & & & & & & & \\ \hline
\end{array}$$
\caption{The frozen variables $a_i$ and cluster variables $x_i$, $z_i$ for $\mathcal A_6$.}
\label{table!e6-orbits}
\end{table}%
We have labelled the Coxeter projection of $\Xi_6$ with these coordinate names, as shown in Figure~\ref{fig!cayley-plane}\,(i) (where the orbit $\{z_1,z_2,z_3\}$ of size three has been squished together in the centre).

{\bf The equations of the Cayley plane.}
The 27 octahedral faces of $\Xi_6$ also split up into two orbits of size 12 and one orbit of size three under the action of the rotation $\sigma$. These three types of octahedral face are shown in Figure~\ref{fig!cayley-plane}\,(ii), and a representative equation from each of these three $\Dih_{12}$-orbits is given by
\begin{align*}
&x_1x_6=a_6x_3 + a_1z_2 +a_8x_4+a_3a_{11}, \tag{a} \\
&x_1x_5= x_3z_3 - a_3x_2 - a_{10}x_4 +a_1a_{12}, \tag{b} \\
&z_1z_2= x_3x_9 + x_6x_{12} + a_2a_8 + a_5a_{11}. \tag{c}
\end{align*}
The choice of $\pm$ sign in front of each monomial in each of these equations is uniquely determined by specifying that all of the equations in orbit (a) obey the analogous positivity rule to equation \eqref{eq!positivity}. Indeed, Macaulay2 agrees that these equations define an irreducible Gorenstein variety $\mathcal{V}_6\subset \PP^{26}$ which has codimension 10 and a Gorenstein resolution with Betti numbers~${( 1, 27, 78, 351, 650, 702, 650, 351, 78, 27, 1 )}$.

\subsection{Finding an initial seed}

We now describe how we found an initial seed for an LPA structure on $\mathcal A_6$. We do this in three steps:
\begin{enumerate}\itemsep=0pt
\item[(1)] First find a candidate for an initial cluster (e.g., an appropriately sized subset of the spinor coordinates for which all other spinor coordinates can be expressed as Laurent polynomials).
\item[(2)] Work out corresponding exchange polynomials for this cluster.
\item[(3)] Check that the corresponding LPA is of finite type and contains all (non-frozen) spinor coordinates on $\mathcal V_6$ as cluster variables.
\end{enumerate}

{\bf Rank of the LPA.}
We first note that we expect such a LPA structure to have rank 5, by Remark~\ref{rmk!seed-rank}.

{\bf Finding an initial cluster.}
We describe three attempts we made in order to find an initial cluster.

{\bf Attempt 1.}
As with the previous cases described in Section~\ref{section:warm-up-cases}, our first thought was to take~${\{x_1,x_2,x_3,x_4,x_5\}}$ as an initial cluster. Unfortunately however, we have the equation
\[x_1x_5 = x_3z_3 - a_3x_2 - a_{10}x_4 + a_1a_{12} \]
that looks like it could make $x_5$ redundant as a cluster variable. It also doesn't help account for the variable $z_3$ appearing in the corresponding exchange polynomial.

{\bf Attempt 2.}
Next we replaced $x_5$ with $z_3$, hoping that $\{x_1,x_2,x_3,x_4,z_3\}$ would work as an initial cluster. Using computer algebra to eliminate variables from the ring $\mathcal A_6$ we discover that all of the other $x_i$ and $z_i$ variables can be written as rational functions in this cluster, but not, unfortunately, as \emph{Laurent polynomials}.

{\bf Attempt 3.}
Finally, although the rational functions obtained in attempt 2 were not Laurent polynomials, a closer inspection reveals that they are \emph{`almost'} Laurent polynomials. Indeed, the denominators are always monomials in the five terms
\[ x_1, x_2, x_3, x_4, x_3z_3 - a_3x_2 - a_{10}x_4. \]
Therefore we introduce $y_3:=x_3z_3 - a_3x_2 - a_{10}x_4$ as a new cluster variable.

\begin{defn}
We let $y_i$ be defined by the expression
\[ y_i := x_iz_i-a_ix_{i-1}-a_{i+7}x_{i+1}, \qquad i\in \ZZ/12\ZZ. \]
\end{defn}

We conclude that we have the following result.

\begin{lem}
We can expand all of the spinors variables $x_i$, $z_i$ as Laurent polynomials in $\{x_1,x_2,x_3,x_4,y_3\}$, and thus we can use this as a candidate for an initial cluster. Moreover, the new variables $y_1,\ldots,y_{12}$ are also all Laurent polynomials in our chosen initial cluster, and they are all distinct. By symmetry we find that both
\[ \{ x_{i-1},x_i,x_{i+1},x_{i+2}, y_i \} \qquad \text{and} \qquad \{ x_{i-1},x_i,x_{i+1},x_{i+2}, y_{i+1} \} \]
are clusters for any value of $i$.
\end{lem}

{\bf Finding the exchange polynomials.}
From equation (b) above we immediately have the relation $x_1x_5 = y_3 + a_1a_{12}$.
Moreover, by manipulations with the equations of types (a), (b), (c), we can write the product $y_2y_3$ as a positive sum of monomials in terms of the frozen coefficients and $x_1$, $x_2$, $x_3$, $x_4$. Indeed we get
\begin{align*}
y_2y_3 &= (x_2z_2 - a_2x_1 - a_9x_3)(x_3z_3 - a_3x_2 - a_{10}x_4) \\
 &= x_1x_4(a_5 x_2 + a_7x_3 + a_2a_{10}) + a_{12}x_3(a_4x_1 + a_1a_9) + a_{12}x_2(a_6x_3 + a_8x_4 + a_3a_{11}).
\end{align*}
Thus we can use these as exchange relations in the following sequence of mutations, together with the expected $\Dih_{12}$ symmetry, to work out what all of the other exchange polynomials should be
\[
\begin{tikzpicture}[scale=0.55]
 \node at (0,1) {$x_1$};
 \node at (1,1) {$x_2$};
 \node at (2,1) {$x_3$};
 \node at (3,1) {$x_4$};
 \node at (1,0) {$y_2$};

 \node at (0+6,1) {$x_1$};
 \node at (1+6,1) {$x_2$};
 \node at (2+6,1) {$x_3$};
 \node at (3+6,1) {$x_4$};
 \node at (2+6,0) {$y_3$};

 \node at (0+12,1) {$x_2$};
 \node at (1+12,1) {$x_3$};
 \node at (2+12,1) {$x_4$};
 \node at (3+12,1) {$x_5$};
 \node at (1+12,0) {$y_3$};

 \node at (0+18,1) {$x_2$};
 \node at (1+18,1) {$x_3$};
 \node at (2+18,1) {$x_4$};
 \node at (3+18,1) {$x_5$};
 \node at (2+18,0) {$y_4$};

 \node at (-3,0.5) {$\cdots$};
 \node at (24,0.5) {$\cdots$.};

 \draw[<->] (4-6 ,0.5) -- (5-6 ,0.5);
 \draw[<->] (4 ,0.5) -- (5 ,0.5);
 \draw[<->] (4+6 ,0.5) -- (5+6 ,0.5);
 \draw[<->] (4+12,0.5) -- (5+12,0.5);
 \draw[<->] (4+18,0.5) -- (5+18,0.5);
 \draw [thick, decorate, decoration = {calligraphic brace}] (-0.5,-0.5) -- (-0.5, 1.5);
 \draw [thick, decorate, decoration = {calligraphic brace}] ( 3.5, 1.5) -- ( 3.5,-0.5);
 \draw [thick, decorate, decoration = {calligraphic brace}] (-0.5+6,-0.5) -- (-0.5+6, 1.5);
 \draw [thick, decorate, decoration = {calligraphic brace}] ( 3.5+6, 1.5) -- ( 3.5+6,-0.5);
 \draw [thick, decorate, decoration = {calligraphic brace}] (-0.5+12,-0.5) -- (-0.5+12, 1.5);
 \draw [thick, decorate, decoration = {calligraphic brace}] ( 3.5+12, 1.5) -- ( 3.5+12,-0.5);
 \draw [thick, decorate, decoration = {calligraphic brace}] (-0.5+18,-0.5) -- (-0.5+18, 1.5);
 \draw [thick, decorate, decoration = {calligraphic brace}] ( 3.5+18, 1.5) -- ( 3.5+18,-0.5);
\end{tikzpicture}
\]
Doing this, we arrive at the following candidate for our the initial seed.

\begin{prop}\label{prop!initial-seed}
The following seed $S$ is an initial seed for a LPA structure on $\mathcal A_6=\CC[\mathcal V_6]$, which is compatible with the $\Dih_{12}$-symmetry:
\[ S = \begin{cases}
x_1, & y_3 + a_{12}a_1, \\
x_2, & a_2x_1(y_3 + a_{10}x_4) + a_9x_3(y_3 + a_1a_{12}) + x_1x_3(a_7x_4 + a_4a_{12}), \\
x_3, & y_3 + a_3x_2 + a_{10}x_4, \\
x_4, & a_{11}(y_3 + a_3x_2) + x_3(a_4x_1 + a_6x_2 + a_1a_9), \\
y_3, & x_1x_4(a_5 x_2 + a_7x_3 + a_2a_{10}) + a_{12}x_3(a_4x_1 + a_1a_9)\\
&\quad + a_{12}x_2(a_6x_3 + a_8x_4 + a_3a_{11}).
\end{cases}\]
\end{prop}

\subsection{Summary of the LPA structure}

{\bf Number of seeds and the exchange graph}
Once we are given the right initial seed it is easy to plug into our code and verify that it generates a LPA of finite type.

\begin{figure}[t]
\centering
\begin{tikzpicture}[xscale=2]
 \node (g) at (0,3) {A};

 \node (d) at (1,4) {B};
 \node (f) at (1,3) {C};
 \node (n) at (1,2) {D};

 \node (e) at (2,5) {E};
 \node (c) at (2,4) {F};
 \node (o) at (2,3) {G};
 \node (b) at (2,2) {H};
 \node (j) at (2,1) {I};

 \node (a) at (3,5) {J};
 \node (m) at (3,4) {K};
 \node (l) at (3,3) {L};
 \node (h) at (3,2) {M};
 \node (k) at (3,1) {N};

 \node (i) at (4,3) {O};

 \draw (a) -- (e) -- (d) -- (g) -- (n) -- (j) -- (b) -- (f) -- (n);
 \draw (g) -- (f) -- (c) -- (d) (c) -- (e) (o) -- (f);
 \draw (c) -- (l) -- (h) (k) -- (i) -- (h) -- (b) -- (k);
 \draw (e) -- (m) -- (l) (m) -- (i) (o) -- (k);

 \draw (e) to [in=60, out=-60] (o);

 \draw (a) to [in=135, out=45,loop] (a);
 \draw (d) to [in=135, out=45,loop] (d);
 \draw (g) to [in=135, out=45,loop] (g);

\end{tikzpicture}
\caption{The $\Dih_{12}$-quotient of the exchange graph of $\mathcal A_6$. Beginning with the initial seed $S$ of Proposition~\ref{prop!initial-seed}, which is in the orbit A, the remaining orbits are named alphabetically, according to the order in which they were found during a breadth-first search of the exchange graph.}\label{fig!exchange-graph}
\end{figure}
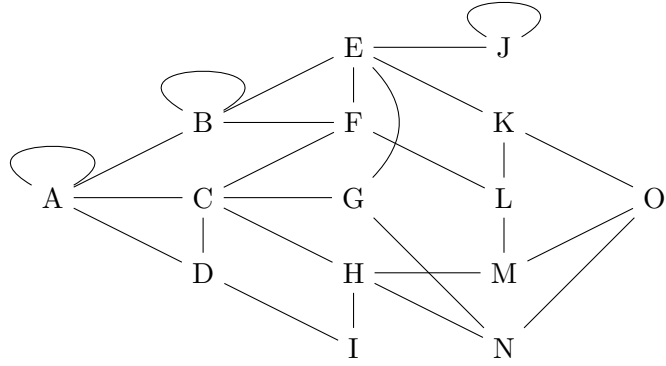

\begin{thm}\label{thm!main-result}
The LPA structure on $\mathcal A_6$, generated by the initial seed of Proposition~{\rm\ref{prop!initial-seed}}, has finite type. In particular it has $264$ seeds and $32$ cluster variables.

The cluster variables consist of the $15$ spinor coordinates $x_1,\ldots,x_{12}$, $z_1$, $z_2$, $z_3$ on $\mathcal V_6$, plus~$17$ additional cluster variables $y_1,\ldots,y_{12}$, $t_1$, $t_2$, $t_3$, $u_1$, $u_2$ where
\begin{enumerate}\itemsep=0pt
\item[$(1)$] $y_1,\ldots,y_{12}$ are quadratics in the spinor variables introduced above,
\item[$(2)$] $t_1$, $t_2$, $t_3$ are quartics in the spinor variables determined by the $\Dih_{12}$-conjugates of the equation
\[ x_1x_4x_7x_{10} = t_1 + a_3a_6a_9a_{12}, \]
\item[$(3)$] $u_1$, $u_2$ are cubics in the spinor variables determined by the $\Dih_{12}$-conjugates of the equation\footnote{It is clear from the equation that $u_1$ is invariant under the shift $i\mapsto i+4$ for $i\in \ZZ/12\ZZ$, but in fact, as a~consequence of the other relations, it turns out that it is invariant under $i\mapsto i+2$ too.}
\[ x_1x_5x_9 = u_1 + a_4a_5x_1 + a_8a_9x_5 + a_{12}a_1x_9 + a_1a_5a_9 + a_4a_8a_{12}. \]
\end{enumerate}
Up to the $\Dih_{12}$-symmetry there are $15$ different orbits of seeds; seven orbits have length $24$, which we label A,$\ldots$,G, and eight orbits have length $12$, which we name H,$\ldots$,O. They are related by mutation according to Table~{\rm\ref{tab!seeds}}, and the $\Dih_{12}$-quotient of the exchange graph is presented in Figure~{\rm\ref{fig!exchange-graph}}.
\end{thm}

\begin{table}[th!]\renewcommand{\arraystretch}{1.2}
\centering
\begin{tabular}{|c|c|} \hline
A & $\begin{array}{ccccc}
x_1&x_2&x_3&x_4&y_2 \\ \hline
y_4&z_2&y_{12}&x_{12}&y_3 \\
\text D&\text B&\text C&\text A&\text A
\end{array}$ \\ \hline
B & $\begin{array}{ccccc}
x_1&x_2&x_4&y_3&z_3 \\ \hline
x_5&x_7&x_{11}&y_{12}&x_3 \\
\text B&\text E&\text F&\text F&\text A
\end{array}$ \\ \hline
C & $\begin{array}{ccccc}
x_1&x_2&x_4&y_{12}&y_2 \\ \hline
u_2&x_{10}&x_{12}&x_3&z_3 \\
\text H&\text G&\text D&\text A&\text F
\end{array}$ \\ \hline
D & $\begin{array}{ccccc}
x_1&x_2&x_3&y_1&y_3 \\ \hline
x_5&u_1&x_{11}&x_{4}&x_{12} \\
\text C&\text I&\text C&\text A&\text A
\end{array}$ \\ \hline
E & $\begin{array}{ccccc}
x_1&x_4&x_7&y_3&z_3 \\ \hline
x_5&y_9&x_2&t_1&y_5 \\
\text F&\text K&\text B&\text J&\text G
\end{array}$ \\ \hline
\end{tabular}
\begin{tabular}{|c|c|} \hline
F & $\begin{array}{ccccc}
x_1&x_2&x_4&y_{12}&z_3 \\ \hline
y_6&x_{10}&x_{11}&y_3&y_2 \\
\text L&\text E&\text B&\text B&\text C
\end{array}$ \\ \hline
G & $\begin{array}{ccccc}
x_1&x_4&x_7&y_3&y_5 \\ \hline
x_5&u_1&x_3&z_2&z_3 \\
\text C&\text N&\text C&\text E&\text E
\end{array}$ \\ \hline
H & $\begin{array}{ccccc}
x_1&x_3&y_3&y_5&u_1 \\ \hline
x_5&x_7&y_{11}&y_1&x_4 \\
\text I&\text N&\text M&\text I&\text C
\end{array}$ \\ \hline
I & $\begin{array}{ccccc}
x_1&x_3&y_1&y_3&u_1 \\ \hline
x_5&x_{11}&y_5&y_{11}&x_2 \\
\text H&\text H&\text H&\text H&\text D
\end{array}$ \\ \hline
J & $\begin{array}{ccccc}
x_1&x_4&x_7&z_2&t_1 \\ \hline
x_{10}&x_{10}&x_{10}&z_3&y_5 \\
\text J&\text J&\text J&\text J&\text E
\end{array}$ \\ \hline
\end{tabular}
\begin{tabular}{|c|c|} \hline
K & $\begin{array}{ccccc}
x_1&x_7&y_5&y_{11}&z_2 \\ \hline
x_9&x_3&x_{10}&x_4&u_2 \\
\text L&\text L&\text E&\text E&\text O
\end{array}$ \\ \hline
L & $\begin{array}{ccccc}
x_1&x_3&y_5&y_{11}&z_2 \\ \hline
x_9&x_7&x_{12}&x_4&u_1 \\
\text K&\text K&\text F&\text F&\text M
\end{array}$ \\ \hline
M & $\begin{array}{ccccc}
x_1&x_3&y_5&y_{11}&u_1 \\ \hline
x_9&x_7&y_1&y_3&z_2 \\
\text O&\text O&\text H&\text H&\text L
\end{array}$ \\ \hline
N & $\begin{array}{ccccc}
x_1&x_7&y_3&y_5&u_1 \\ \hline
x_5&x_3&y_{11}&y_9&x_4 \\
\text H&\text H&\text O&\text O&\text G
\end{array}$ \\ \hline
O & $\begin{array}{ccccc}
x_1&x_7&y_5&y_{11}&u_1 \\ \hline
x_9&x_3&y_9&y_3&z_2 \\
\text M&\text M&\text N&\text N&\text K
\end{array}$ \\ \hline
\end{tabular}
\caption{Representatives for each of the 15 orbits of seeds A,$\ldots$,O. The top row of each entry contains the five cluster variables in the seed. The second row records which cluster variable is obtained by mutating the seed at the variable directly above it, leading to a seed in the orbit given by the label on the third row.}\label{tab!seeds}
\end{table}

We can also check various things about the structure of the exchange graph, such as the fact that every $x$ (resp.\ $y$, $z$, $t$, $u$) variable belongs to 60 (resp.\ 32, 40, 8, 36) seeds.

{\bf Positivity.}
By inspecting the output of our computation, which consists of all of the seeds for $\mathcal A_6$ (including the Laurent expansion of all of the cluster variables), we have the following result.

\begin{cor}\label{cor!positivity}
The positivity phenomenon holds for $\mathcal A_6$. In other words, all of the coefficients in the Laurent expansions of the cluster variables and all the coefficients in exchange polynomials of each seed are positive.
\end{cor}

This is somewhat unexpected, since enforcing the positivity rule in equation \eqref{eq!positivity} on the equations of type (a) necessarily creates minus signs in some of the other spinor equations defining the Cayley plane, e.g.,
\[ x_1x_5 = x_3z_3 - a_3x_2 - a_{10}x_4 + a_1a_{12}. \]
However, in order to get a Laurent phenomenon we needed to introduce the new cluster variable~${y_3=x_3z_3 - a_3x_2 - a_{10}x_4}$ and this allows us to rewrite this equation with positive coefficients as $x_1x_5 = y_3 + a_1a_{12}$.

\begin{rmk}
Positivity was proved for cluster algebras by Gross, Hacking, Keel and Kontsevich~\cite{ghkk} by associating a \emph{consistent scattering diagram} to a cluster algebra. The proof follows by interpreting the coefficients in the Laurent expansion of each cluster monomial as a count of broken lines in the scattering diagram. The consistent scattering diagram for the LPA structure on $\mathcal A_5$ was constructed in \cite{duc}, and it should be possible to construct one for $\mathcal A_6$ in a similar manner.
\end{rmk}

{\bf A final remark.}
We recall that the definition of mutation in an LPA uses the exchange Laurent polynomials $\widehat{F}_i$, rather than the exchange polynomials $F_i$. We always have $F_i=\widehat F_i$ in the case of cluster algebras, and it is tempting to think we might be able to dispense with the~$\widehat F_i$ in the general case of an LPA. However, it is crucial to work with the $\widehat F_i$ in order for the LPA structure on $\mathcal A_6$ to have finite type. Moreover, this LPA provides an example for which every seed has at least one direction $i$ in which $F_i\neq\widehat F_i$.

\section{The Freudenthal variety}

We conjecture the existence of a similar LPA structure on the homogeneous coordinate ring of the Freudenthal variety $\mathcal V_7$.

\begin{conj}
There is a finite type LPA structure of rank $10$ on $\mathcal A_7=\CC[\mathcal V_7]$, which has the positivity phenomenon.
\end{conj}

The LPA should have rank 10 by Remark~\ref{rmk!seed-rank}. We get as far as writing down the equations for $\mathcal V_7$ (as we did for $\mathcal V_6$ in Section~\ref{sec!cayley-eqns}). 

\subsection{The equations of the Freudenthal variety}
The Freudenthal variety has an embedding $\mathcal V_7\!\!\subset\!\! \PP^{55}$ where the 56 variables $a_1,\ldots,a_{18}$, $x_1,\ldots,x_{18}$, $y_1,\ldots,y_{18}$, $z_1$, $z_2$ can be put into one-to-one correspondence with the 56 lines on a del Pezzo surface of degree 2. The Coxeter rotation has order $18$ and the action on the variables is on the labels.

There are seven orbits of equations:
\begin{gather}
x_1x_2 		= a_3x_{17} + a_2y_0 + a_1y_3 + a_0x_4 + a_{17}a_4, 				\tag{a}\\
y_{14}y_1 	= x_{15}y_4 - a_0x_{10} - a_{15}x_5 + x_0y_{11} + x_{13}x_2, 	\tag{b} \\
y_2y_1 		= a_2y_{12} + y_5x_0 + y_{16}x_3 + a_1y_9 + a_5a_{16} ,			\tag{c} \\
x_4y_1 		= x_2y_5 + a_5x_{17} + a_4y_{16} - a_2y_8 - a_1x_7 ,				\tag{d} \\
x_2x_6 		= x_4y_4 \mp a_4z_2 - a_5y_0 - a_3y_8 + a_1a_7 ,					\tag{e} \\
y_1y_3 		= x_{17}x_5 \pm x_2z_1 + a_0x_7 + a_4x_{15} - a_2y_{11},			\tag{f} \\
y_{4}y_{16} = \pm y_1z_2 - y_8x_0 - y_{12}x_2 - a_1x_{10} + a_5a_{15}, 			\tag{g}
\end{gather}
corresponding to the diagrams
\begin{center}\begin{tikzpicture}[scale=1.5,rotate=40]
 \setupnodes
 \node at ({1.3*cos(135-40)},{1.3*sin(135-40)}) {(a)};
 \node[red] at (a0) {$\bullet$};
 \node[red] at (a5) {$\bullet$};
 \node[blue] at (a1) {$\bullet$};
 \node[blue] at (b5) {$\bullet$};
 \node[green] at (a2) {$\bullet$};
 \node[green] at (c4) {$\bullet$};
 \node[orange] at (a3) {$\bullet$};
 \node[orange] at (c1) {$\bullet$};
 \node[purple] at (a4) {$\bullet$};
 \node[purple] at (b0) {$\bullet$};
 \node[black] at (b2) {$\bullet$};
 \node[black] at (b3) {$\bullet$};
\end{tikzpicture}$\quad$\begin{tikzpicture}[scale=1.5,rotate=20]
 \setupnodes
 \node at ({1.3*cos(135-20)},{1.3*sin(135-20)}) {(b)};
 \node[red] at (c1) {$\bullet$};
 \node[red] at (c6) {$\bullet$};
 \node[blue] at (b5) {$\bullet$};
 \node[blue] at (c16) {$\bullet$};
 \node[green] at (a2) {$\bullet$};
 \node[green] at (b10) {$\bullet$};
 \node[orange] at (a5) {$\bullet$};
 \node[orange] at (b15) {$\bullet$};
 \node[purple] at (b2) {$\bullet$};
 \node[purple] at (c9) {$\bullet$};
 \node[black] at (b0) {$\bullet$};
 \node[black] at (b7) {$\bullet$};
\end{tikzpicture}$\quad$\begin{tikzpicture}[scale=1.5]
 \setupnodes
 \node at ({1.3*cos(135)},{1.3*sin(135)}) {(c)};
 \node[red] at (b3) {$\bullet$};
 \node[red] at (c8) {$\bullet$};
 \node[blue] at (a4) {$\bullet$};
 \node[blue] at (c12) {$\bullet$};
 \node[green] at (a5) {$\bullet$};
 \node[green] at (c15) {$\bullet$};
 \node[orange] at (a1) {$\bullet$};
 \node[orange] at (a8) {$\bullet$};
 \node[purple] at (c4) {$\bullet$};
 \node[purple] at (c5) {$\bullet$};
 \node[black] at (c1) {$\bullet$};
 \node[black] at (b6) {$\bullet$};
\end{tikzpicture}$\quad$\begin{tikzpicture}[scale=1.5,rotate=-60]
 \setupnodes
 \node at ({1.3*cos(135+60)},{1.3*sin(135+60)}) {(d)};
 \node[red] at (a2) {$\bullet$};
 \node[red] at (c8) {$\bullet$};
 \node[blue] at (a4) {$\bullet$};
 \node[blue] at (c16) {$\bullet$};
 \node[green] at (a5) {$\bullet$};
 \node[green] at (b17) {$\bullet$};
 \node[orange] at (a1) {$\bullet$};
 \node[orange] at (b7) {$\bullet$};
 \node[purple] at (b2) {$\bullet$};
 \node[purple] at (c5) {$\bullet$};
 \node[black] at (c1) {$\bullet$};
 \node[black] at (b4) {$\bullet$};
\end{tikzpicture}

\begin{tikzpicture}[scale=1.5,rotate=-80]
 \setupnodes
 \node at ({1.3*cos(135+80)},{1.3*sin(135+80)}) {(e)};
 \node[red] at (a3) {$\bullet$};
 \node[red] at (c8) {$\bullet$};
 \node[blue] at (a4) {$\bullet$};
 \node[blue] at (0,0) {$\bullet$};
 \node[green] at (a5) {$\bullet$};
 \node[green] at (c0) {$\bullet$};
 \node[orange] at (a1) {$\bullet$};
 \node[orange] at (a7) {$\bullet$};
 \node[purple] at (b2) {$\bullet$};
 \node[purple] at (b6) {$\bullet$};
 \node[black] at (c4) {$\bullet$};
 \node[black] at (b4) {$\bullet$};
\end{tikzpicture}$\quad$\begin{tikzpicture}[scale=1.5,rotate=-60]
 \setupnodes
 \node at ({1.3*cos(135+60)},{1.3*sin(135+60)}) {(f)};
 \node[red] at (b3) {$\bullet$};
 \node[red] at (0,0) {$\bullet$};
 \node[blue] at (a3) {$\bullet$};
 \node[blue] at (c12) {$\bullet$};
 \node[green] at (a5) {$\bullet$};
 \node[green] at (b16) {$\bullet$};
 \node[orange] at (a1) {$\bullet$};
 \node[orange] at (b8) {$\bullet$};
 \node[purple] at (c4) {$\bullet$};
 \node[purple] at (c2) {$\bullet$};
 \node[black] at (b0) {$\bullet$};
 \node[black] at (b6) {$\bullet$};
\end{tikzpicture}$\quad$\begin{tikzpicture}[scale=1.5,rotate=-120]
 \setupnodes
 \node at ({1.3*cos(135+120)},{1.3*sin(135+120)}) {(g)};
 \node[red] at (a6) {$\bullet$};
 \node[red] at (b15) {$\bullet$};
 \node[blue] at (b5) {$\bullet$};
 \node[blue] at (c13) {$\bullet$};
 \node[green] at (a2) {$\bullet$};
 \node[green] at (a10) {$\bullet$};
 \node[orange] at (b7) {$\bullet$};
 \node[orange] at (c17) {$\bullet$};
 \node[purple] at (c6) {$\bullet$};
 \node[purple] at (0,0) {$\bullet$};
 \node[black] at (c3) {$\bullet$};
 \node[black] at (c9) {$\bullet$};
\end{tikzpicture}\end{center}
and some additional quadratic equations
\begin{align*} z_1z_2 ={}& y_3y_{12} + y_7y_{16} + y_8y_{17} + x_1x_{10} + x_3x_{12} + x_5x_{14} + a_0a_9 + a_2a_{11} - a_3a_{12}\\
& + a_4a_{13} + a_6a_{15} \end{align*}
that are implied by the others. Here the choice of $\pm$ sign in the equations is due to the fact that the $z_i$ variables do not appear in the equations of type (a), and thus the positivity rule of equation \eqref{eq!positivity} does not determine the signs in front of the monomials containing exactly one~$z_i$.

From our observations on the previous cases, it seems like it will be constructive to introduce new cluster variables which will allow us to rearrange the equations so that they satisfy the positivity phenomenon, such as
\begin{align*}
&t_?:= y_{14}y_1 - x_0y_{11} - x_{13}x_2 = x_{15}y_4 - a_0x_{10} - a_{15}x_5, \\
&u_?:= x_4y_1 - a_5x_{17} -a_4y_{16} = x_2y_5 - a_2y_8 - a_1x_7,
\end{align*}
and so on.

However, at this point we get stuck. Trying to follow our previous approach of identifying an initial cluster, by using computer algebra to eliminate variables in a ring of codimension 28, proves to be a step too far.

\appendix

\section{Sage code} \label{appendixA}
The sequence of LPAs associated to the spaces in this paper can be computed using the Sage package \textsc{LPASeed} \cite{sagecode}. To define the initial seed for the $E_4$ case, one declares the variable names and corresponding polynomials as follows:
\begin{lstlisting}
sage: var('x1, x2')
(x1, x2)
sage: var('a2, a4, a5, a1, a3')
(a2, a4, a5, a1, a3)
sage: S = LPASeed({x1: a2*x2 + a4*a5, x2: a1*x1 + a3*a4},
....: coefficients=[a2, a4, a5, a1, a3])
\end{lstlisting}

Then one can analyse the structure of the mutations using the built-in methods. For instance, to get the number of seeds in the exchange graph:
\begin{lstlisting}
sage: len(S.mutation_class())
5
\end{lstlisting}

A list of all the cluster variables:
\begin{lstlisting}
sage: S.variable_class()
[(x1*a5*a1 + x2*a2*a3 + a4*a5*a3)/(x1*x2),
 (x1*a1 + a4*a3)/x2,
 (x2*a2 + a4*a5)/x1,
 x2,
 x1]
\end{lstlisting}

We can work out larger examples rather quickly and even plot their exchange graphs. The LPA for the orthogonal $\text{OGr}(5, 10)$ case corresponding to $E_5$ can be fully understood in a few lines:
\begin{lstlisting}
sage: var('x1, x2, x3')
(x1, x2, x3)
sage: var('a1, a2, a3, a4, a5, a6, a7, a8')
(a1, a2, a3, a4, a5, a6, a7, a8)
sage: F1 = a5*x2 + a8*x3 + a2*a3
sage: F2 = a6*x1*x3 + a3*a4*x1 + a8*a1*x3 + a1*a2*a3
sage: F3 = a4*x1 + a7*x2 + a1*a2
sage: S = LPASeed({x1: F1, x2: F2, x3: F3},
....: coefficients=[a1,a2,a3,a4,a5,a6,a7,a8])
sage: len(S.mutation_class())
16
sage: len(S.variable_class())
10
sage: show(S.exchange_graph())
\end{lstlisting}
The final command shows the exchange graph featured in Figure~\ref{fig:exchange-graph} to the user.

\subsection*{Acknowledgements}

Both authors would like to thank Anna Felikson for some very useful discussions whilst undertaking this research project. We are also grateful to the anonymous referees for their detailed feedback.

\pdfbookmark[1]{References}{ref}
\LastPageEnding


\begin{thebibliography}{99}
\footnotesize\itemsep=0pt

\bibitem{c}
Coxeter H.S.M., Regular polytopes, 3rd ed., Dover Publications, Inc., New York,
 1973.

\bibitem{sagecode}
Daisey O., Laurent phenomenon algebra seed cell, available online at
 \url{https://oliverdaisey.github.io/code.html}.

\bibitem{dol}
Dolgachev I.V., Classical algebraic geometry. {A} modern view, \href{https://doi.org/10.1017/CBO9781139084437}{Cambridge
 University Press}, Cambridge, 2012.

\bibitem{duc}
Ducat T., The 3-dimensional {L}yness map and a self-mirror log {C}alabi--{Y}au
 3-fold, \href{https://doi.org/10.1007/s00229-023-01497-0}{\textit{Manuscripta Math.}} \textbf{174} (2024), 87--140,
 \href{https://arxiv.org/abs/2105.07843}{arXiv:2105.07843}.

\bibitem{GLS}
Gei{\ss} C., Leclerc B., Schr\"oer J., Partial flag varieties and preprojective
 algebras, \href{https://doi.org/10.5802/aif.2371}{\textit{Ann. Inst. Fourier (Grenoble)}} \textbf{58} (2008),
 825--876, \href{https://arxiv.org/abs/math.RT/0609138}{arXiv:math.RT/0609138}.

\bibitem{ghk}
Gross M., Hacking P., Keel S., Birational geometry of cluster algebras,
 \href{https://doi.org/10.14231/AG-2015-007}{\textit{Algebr. Geom.}} \textbf{2} (2015), 137--175, \href{https://arxiv.org/abs/1309.2573}{arXiv:1309.2573}.

\bibitem{ghkk}
Gross M., Hacking P., Keel S., Kontsevich M., Canonical bases for cluster
 algebras, \href{https://doi.org/10.1090/jams/890}{\textit{J.~Amer. Math. Soc.}} \textbf{31} (2018), 497--608,
 \href{https://arxiv.org/abs/1411.1394}{arXiv:1411.1394}.

\bibitem{lp}
Lam T., Pylyavskyy P., Laurent phenomenon algebras, \href{https://doi.org/10.4310/CJM.2016.v4.n1.a2}{\textit{Camb.~J. Math.}}
 \textbf{4} (2016), 121--162, \href{https://arxiv.org/abs/1206.2611}{arXiv:1206.2611}.

\bibitem{m}
Moufang R., Alternativk\"orper und der {S}atz vom vollst\"andigen {V}ierseit
 {$(D_9)$}, \href{https://doi.org/10.1007/BF02940648}{\textit{Abh. Math. Sem. Univ. Hamburg}} \textbf{9} (1933),
 207--222.

\bibitem{scott}
Scott J.S., Grassmannians and cluster algebras, \href{https://doi.org/10.1112/S0024611505015571}{\textit{Proc. London Math.
 Soc.~(3)}} \textbf{92} (2006), 345--380, \href{https://arxiv.org/abs/math/0311148}{arXiv:math/0311148}.

\bibitem{snow}
Snow D.M., Homogeneous vector bundles, in Group {A}ctions and {I}nvariant
 {T}heory, \textit{CMS Conf. Proc.}, Vol.~10, American Mathematical Society,
 Providence, RI, 1989, 193--205.

\bibitem{Spacek-Wang}
Spacek P., Wang C., Towards {L}andau--{G}inzburg models for cominuscule spaces
 via the exceptional cominuscule family, \href{https://doi.org/10.1016/j.jalgebra.2023.03.039}{\textit{J.~Algebra}} \textbf{630}
 (2023), 334--393, \href{https://arxiv.org/abs/2204.03548}{arXiv:2204.03548}.

\end{thebibliography}
\end{document}